\documentclass{amsart}

\usepackage{amsfonts,amssymb,verbatim,amsmath,amsthm,latexsym,textcomp,amscd}
\usepackage{latexsym,amsfonts,amssymb,epsfig,verbatim}
\usepackage{amsmath,amsthm,amssymb,latexsym,graphics,textcomp}
\usepackage{graphicx}
\usepackage{color}
\usepackage{url}

\input xy
\xyoption{all}

\theoremstyle{plain}
\newtheorem{theorem}{Theorem}[section]

\newtheorem{question}[theorem]{Question}

\newtheorem{scholium}[theorem]{Scholium}
\newtheorem{obs}[theorem]{Observation}

\newtheorem{prop}[theorem]{Proposition}
\newtheorem{lemma}[theorem]{Lemma}
\newtheorem{cor}[theorem]{Corollary}
\newtheorem{defn}[theorem]{Definition}

\newtheorem{rmk}[theorem]{Remark}

\newcommand{\bdy}{\partial}
\newcommand{\C}{{\mathbb C}}

\newcommand{\natls}{{\mathbb N}}

\newcommand{\reals}{{\mathbb R}}

\newcommand\Hyp{{\mathbf H}}

\newcommand\Z{{\mathbb Z}}
\newcommand\R{{\mathbb R}}

\newcommand\AAA{{\mathcal A}}
\newcommand\BB{{\mathcal B}}

\newcommand\EE{{\mathcal E}}
\newcommand\FF{{\mathcal F}}
\newcommand\GG{{\mathcal G}}
\newcommand\HH{{\mathcal H}}

\newcommand\JJ{{\mathcal J}}
\newcommand\KK{{\mathcal K}}
\newcommand\LL{{\mathcal L}}
\newcommand\MM{{\mathcal M}}

\newcommand\OO{{\mathcal O}}
\newcommand\PP{{\mathcal P}}
\newcommand\QQ{{\mathcal Q}}

\newcommand\VV{{\mathcal V}}

\newcommand\PMF{{\PP\kern-2pt\MM\FF}}

\newcommand\PML{{\PP\kern-2pt\MM\LL}}

\newcommand{\fsubd}{\mathrel{{\scriptstyle\searrow}\kern-1ex^d\kern0.5ex}}
\newcommand{\bsubd}{\mathrel{{\scriptstyle\swarrow}\kern-1.6ex^d\kern0.8ex}}
\newcommand{\fsubeq}{\mathrel{\raise-.7ex\Hbox{$\overset{\searrow}{=}$}}}
\newcommand{\bsubeq}{\mathrel{\raise-.7ex\Hbox{$\overset{\swarrow}{=}$}}}

\begin{document}

\title{Pattern Rigidity and the Hilbert-Smith Conjecture}
\author{Mahan Mj}

\thanks{Research partly supported by  CEFIPRA research project grant 4301-1}

\begin{abstract} 
In this paper we initiate a study of the topological group $PPQI(G,H)$ of pattern-preserving
quasi-isometries for $G$  a hyperbolic Poincare duality group and $H$ an  infinite quasiconvex subgroup
of infinite index in $G$. Suppose  
$\partial G$ admits a visual metric $d$ with 
$dim_{haus} < dim_t +2$, where $dim_{haus}$ is the
Hausdorff dimension and $dim_t$ is the
topological dimension of $(\partial G,d)$. Equivalently suppose that $ACD(\partial G) < dim_t + 2$ where 
$ACD(\partial G)$ denotes the Ahlfors regular conformal dimension of $\partial G$.\\
a) If $Q_u$ is a group 
of pattern-preserving uniform quasi-isometries (or more generally any locally compact group
of pattern-preserving  quasi-isometries) containing $G$, then $G$ is of finite index in $Q_u$. \\
b) If instead, $H$ is a codimension one filling subgroup, and $Q$ is any group 
of  pattern-preserving  quasi-isometries containing $G$, then $G$ is of finite index in $Q$. Moreover,
 {\bf (Topological Pattern Rigidity)}
if $L$ is the limit set of $H$, $\LL$ is the collection of translates of $L$ under $G$, 
and  $Q$ is any pattern-preserving group
of {\it homeomorphisms} of $\partial G$  preserving $\LL$ and containing $G$, 
then  the index of $G$ in $Q$ is finite.\\ 
We find analogous results in the realm of relative hyperbolicity, regarding an equivariant
collection of horoballs as a symmetric pattern in the universal cover of a complete finite volume non-compact manifold of
pinched negative curvature.
Combining our main result with a theorem of Mosher-Sageev-Whyte, we obtain QI rigidity results.

An important ingredient of the proof is a version of the Hilbert-Smith conjecture for certain metric
measure spaces, which uses the full strength of Yang's theorem on actions of the p-adic integers on
homology manifolds. This might be of independent interest.

\smallskip

\begin{center}

{\em AMS subject classification =   20F67(Primary), 22E40   57M50}

\end{center}

\end{abstract}

\maketitle

\tableofcontents

\section{Preliminaries}

\subsection{Statement of Results} 
In this paper we start studying the full group of `pattern-preserving quasi-isometries' for pairs
$(G,H)$ where $G$ is a (Gromov) hyperbolic group and $H$ an infinite quasiconvex subgroup of
infinite index in $G$.
In \cite{gromov-ai} Gromov proposed the project of classifying
 finitely generated groups up to quasi-isometry, as well as the study of the group $QI(X)$
of quasi-isometries of a space $X$, where two quasi-isometries  are identified if they
lie at a bounded distance from each other. A class of
groups where
any two members are quasi-isometric if and only if they are commensurable
is said to be quasi-isometrically rigid.
However, any
class
of groups acting freely, cocompactly and properly discontinuously on some fixed proper hyperbolic 
metric space $\Hyp$ are
quasi-isometric to $\Hyp$ and hence
to each other. In this context (or in a context where quasi-isometric
rigidity is not known) it makes sense to ask a relative version
of Gromov's question. To obtain rigidity results, we 
 impose additional restrictions on the quasi-isometries by requiring that they preserve some 
additional structure given by a `symmetric pattern' of subsets. A 
`symmetric pattern' of subsets
 roughly means a $G$-equivariant  collection $\JJ$ of convex (or uniformly quasiconvex)
cocompact subsets
in $\Hyp$ (see Section \ref{pats} for detailed definitions). Then the relative version
of Gromov's question for classes of pairs $(G,H)$ was formulated by Mosher-Sageev-Whyte \cite{msw2}
as the following {\it pattern rigidity} question: \\

\begin{question} Given a quasi-isometry $q$ of two such pairs
$(G_i, H_i)$ ($i = 1, 2$) pairing a $(G_1, H_1)$-symmetric pattern
$\JJ_1$ with a $(G_2, H_2)$-symmetric pattern
$\JJ_2$, does there exist an abstract commensurator $I$ which performs the
same pairing? 
\label{motvn}
\end{question}

The study of this question was initiated  by Schwartz \cite{schwartz_pihes},  \cite{schwarz-inv}, where $G$
is a lattice in a rank one symmetric space. 
The paper \cite{schwartz_pihes} deals with 
symmetric
patterns of convex sets (horoballs) whose limit sets are single points, and 
\cite{schwarz-inv} deals with 
symmetric
patterns of convex sets (geodesics)
whose limit sets consist of two points.  In \cite{biswas-mj}, Biswas and Mj generalized Schwartz' result
to certain Duality and PD subgroups of rank one symmetric spaces. In \cite{bi}, Biswas completely
solved the pattern rigidity problem for $G$ a uniform lattice in real hyperbolic space
and $H$ any infinite quasiconvex subgroup of infinite index in $G$. However, all these papers used, in an
essential way, the linear structure of the groups involved, and the techniques fail for $G$ the fundamental group
of a general closed negatively curved manifold. (This point is specifically mentioned by Schwartz in
\cite{schwarz-inv}). Further, the study in \cite{schwartz_pihes},  \cite{schwarz-inv}, \cite{biswas-mj},
\cite{bi} boils down to the study of a single pattern-preserving quasi-isometry between pairs
$(G_1,H_1)$ and $(G_2,H_2)$. We propose a different perspective in this paper by studying the full group
$PPQI(G,H)$ of pattern-preserving (self) quasi-isometries of a pair $(G,H)$ for $G$ a hyperbolic group and 
$H$ any infinite quasiconvex subgroup of infinite index. The features of $G$ that we shall use are general enough to
go beyond the linear context while at the same time being strong enough to ensure rigidity in certain contexts. Some of the
ingredients of this paper are: \\
1) The boundary of a Poincare duality (PD for short) hyperbolic group  is a homology manifold 
(cf. Definition \ref{d:gm}) by a Theorem of Bestvina-Mess
\cite{bestvina-mess}.\\
2) The algebraic topology of homology manifolds  imposes restrictions on what kinds of groups may act on them
by Theorems of Newman \cite{newman}, Smith \cite{smith3} and Yang \cite{yang-p}. \\
3) Boundaries of hyperbolic groups equipped with the visual metric also have a metric measure space structure
with the property that they are Ahlfors regular (cf. Definition \ref{ahlforsreg-def}).\\
4) Quasiconformal analysis can be conducted in the general context of Ahlfors regular metric spaces.\\
5) A combinatorial cross-ratio can be constructed on the boundary of a hyperbolic group in the presence
of a codimension one subgroup. (Roughly speaking these are subgroups whose Cayley graphs coarsely separate the
 Cayley graph of the group. See first paragraph of Section \ref{codimone}
for a formal definition of codimension one subgroups.)\\
We refer the reader to \cite{hock-young}, p. 145,
 for details on topological dimension and  \cite{davis-pd}
for details on PD groups.

Of these ingredients, the first two come from (a somewhat forgotten chapter of) algebraic topology, the next two
from a very active new area of analysis on metric measure spaces, while the last comes from geometric
group theory proper. Topological actions of finite groups on manifolds and homological consequences
of actions of $p$-adics on  manifolds form the two main ingredients for a proof of the Hilbert-Smith conjecture
for bi-Lipschitz \cite{repovs} and quasiconformal \cite{martin-hs} actions. We first generalize the result
of Martin \cite{martin-hs} to Ahlfors regular metric spaces that are boundaries of PD hyperbolic groups
and obtain the following.

\smallskip

{\bf Theorem \ref{hs} and Corollary \ref{hscor} :} 
Let $G$ be a Poincare duality hyperbolic group and $Q$ be a group of (boundary values of)
quasi-isometries of $G$. Suppose $d$ is a visual metric on $\partial G$ with 
$ dim_{haus} < dim_t +2$, where $dim_{haus}$ is the
Hausdorff dimension (cf. Definition \ref{haus-def}) and $dim_t$ is the
topological dimension of $(\partial G,d)$. 
Equivalently suppose that $ACD(\partial G) < dim_t + 2$ where $ACD(\partial G)$ denotes the Ahlfors regular conformal dimension 
(cf. Definition \ref{acd}) of $\partial G$. 
Then $Q$ cannot contain a copy of 
 $Z_{(p)}$, where 
$Z_{(p)}$ denotes the $p$-adic integers. Hence if $Q$
 is finite dimensional locally compact, it must be a Lie group.

\smallskip

Theorem \ref{hs} and Corollary \ref{hscor} give a strong affirmative answer to a question of Iwaniec and Martin (\cite{iwaniec-martin}
Remark 2, p. 527).

\smallskip

As in \cite{martin-hs}, there is
no assumption on the uniformity of the quasiconformal maps in $Z_{(p)}$. The analogue of Theorem \ref{hs}
is {\it false} for purely topological actions \cite{raymond-williams} on homology manifolds. Hence the 
quasi-isometry assumption is crucial here. The statement that a topological manifold does not admit
an effective topological $Z_{(p)}$ action is the famous Hilbert-Smith conjecture (but does not imply
Corollary \ref{hscor}).

Corollary \ref{hscor}
will be a crucial ingredient in our approach to pattern rigidity. 
Another property we shall be investigating
in some detail is the notion of `topological infinite divisibility' (see section \ref{infdivs} for definitions).
 The notion we introduce is somewhat
weaker
than related existing notions in the literature. In this generality, we prove

\smallskip

{\bf Propositions \ref{td} and \ref{pt} } Let $G$ be a hyperbolic group and $H$ an  infinite quasiconvex subgroup
of infinite index in $G$. Then \\
a) any group 
of  pattern-preserving  quasi-isometries  is totally disconnected
and contains no topologically infinitely divisible elements. \\
b) If $G$ is a Poincare duality group, the group $QI(G)$ of quasi-isometries cannot contain arbitrarily
small torsion elements. 

\smallskip

We obtain stronger results under the assumption that $G$ is a Poincare duality group (e.g. 
the fundamental group
of a  closed negatively curved manifold) with some restrictions on the visual metric on its boundary.
(We refer the reader to the first paragraph of Section \ref{codimone}
for the notion of codimension one filling subgroups.)\\
{\bf Theorems \ref{weakrig}, \ref{strongrig} and \ref{toprig}} Let $G$ be a hyperbolic Poincare duality group and $H$ an  infinite quasiconvex subgroup
of infinite index in $G$. Suppose further that for some visual metric on $\partial G$, $dim_{haus}(\partial G)
< dim_t(\partial G) +2$, where $dim_{haus}$ and $dim_t$ denote Hausdorff and topological dimension respectively.
Equivalently suppose that $ACD( \partial G) < dim_t + 2$ where $ACD(\partial G)$ denotes the Ahlfors regular conformal dimension of $\partial G$.\\
a) If $Q_u$ is a group 
of pattern-preserving uniform quasi-isometries (or more generally any  locally compact group
of pattern-preserving  quasi-isometries)
containing $G$, then $G$ is of finite index in $Q_u$. \\
b) If further, $H$ is a codimension one filling subgroup, and $Q$ is any group 
of (not necessarily uniform) pattern-preserving  quasi-isometries containing $G$, then $G$ is of finite index in $Q$. \\
c) {\bf Topological Pattern Rigidity} Under the assumptions of (b),
let $L$ be the limit set of $H$ and $\LL$ be the collection of translates of $L$ under $G$.  
Let  $Q$ be any pattern-preserving group
of homeomorphisms of $\partial G$  preserving $\LL$ and containing $G$. Then  the index of $G$ in $Q$ is finite.

\smallskip

Theorem \ref{toprig} is a generalization of a Theorem of Casson-Bleiler \cite{CB} and Kapovich-Kleiner
\cite{kap-kl-lowd} to all dimensions. Casson-Bleiler \cite{CB} and Kapovich-Kleiner
\cite{kap-kl-lowd} proved Theorem \ref{toprig} for $G$ the fundamental group of a surface
and $H$ an infinite cyclic subgroup corresponding to a filling curve.

\smallskip

\noindent {\bf Codimension one filling subgroups:} (See first paragraph of Section \ref{codimone}
for definitions.) The existence of a codimension one filling  quasiconvex
subgroup $H$ (or more generally a finite family of codimension one quasiconvex
subgroups that are filling as a collection) of
 a hyperbolic group $G$ ensures that $G$ acts properly, cocompactly on a CAT(0) cube complex \cite{sageev}
by a recent result of Bergeron-Wise \cite{berg-wise}. Thus, Theorems \ref{weakrig}, \ref{strongrig} and \ref{toprig} deal
with the pattern-rigidity for PD(n) hyperbolic groups acting properly and cocompactly on a CAT(0) cube complex.

\smallskip

We also derive QI rigidity results for fundamental groups
of certain non-compact
negatively curved manifolds of finite volume, by deriving analogues of Theorem \ref{weakrig} for symmetric patterns
of horoballs and combining it with a Theorem of Behrstock-Drutu-Mosher \cite{BDM}. The hypotheses in
the following Theorem are satisfied by fundamental groups of finite volume
complete non-compact manifolds of sufficiently pinched negative curvature and dimension bigger than 2.

\smallskip

{\bf Theorem \ref{weakrigrh}} Let $M = M^n$ be a complete finite volume  manifold of pinched negative curvature with $n > 2$. Let $G= \pi_1(M)$. Suppose that there exists
 a visual metric $d$ on 
$\partial (\widetilde{M})$ with 
$dim_{haus} < dim_t +2$, where $dim_{haus}$ is the
Hausdorff dimension and $dim_t$ is the
topological dimension of $(\partial (\widetilde{M}),d)$.  Equivalently suppose that $ACD(\partial (\widetilde{M})) < dim_t + 2$ where $ACD(\partial
(\widetilde{M}))$ denotes the Ahlfors regular conformal dimension of $\partial (\widetilde{M})$.\\
Let $\Gamma$ be a Cayley graph of $G$ with respect to a finite generating set. Let $Q$ be a group 
of  uniform quasi-isometries of $\Gamma$ containing $G$. Then $G$ is of finite index in $Q$. In particular,
$Q \subset Comm (G)$, where $Comm (G)$ denotes the abstract commensurator of $G$.

\smallskip

The author learnt the following Scholium from Misha Gromov \cite{misha}.

\begin{scholium} \label{gromov}
If two discrete groups can be embedded in the same locally compact group nicely, they are
as good as
commensurable.
\end{scholium}

A partial aim of this paper is to make Scholium \ref{gromov} precise in the context of pattern rigidity. It
follows from Theorems
\ref{weakrig}  and \ref{weakrigrh} that in the context of pattern rigidity or QI rigidity of (fundamental
groups of) finite volume
complete non-compact manifolds of pinched negative curvature and dimension bigger than 2,
`as good as' can be replaced by `actually' in Scholium \ref{gromov}. Thus, Theorems
\ref{weakrig} and \ref{weakrigrh}
reduce the problem of pattern rigidity to the weaker problem of embedding two groups
simultaneously in the same locally compact group and Theorem \ref{strongrig} carries out this embedding
under certain hypotheses.

\bigskip

\noindent {\bf Notation:} \\
To prevent confusion we fix two pieces of notation: \\
1) $\Z_p$ will denote the integers mod $p$. \\
2) $Z_{(p)}$ will denote the $p$-adic integers.

\subsection{Dotted geodesic metric spaces}
\begin{defn} A dotted metric space is a metric space $X$, where $d(x, y)$ is an integer for all $x, y \in X$.
A dotted geodesic metric space is a dotted metric space $X$, such that for all $x, y \in X$, there exists
an isometric map $ \sigma \colon [0, d(x, y)] \cap \Z \rightarrow X$ with $\sigma (0) = x$ and $ \sigma ( d(x, y)) = y$.
A dotted metric space is proper if every ball $N_k(x)$ is finite.
\end{defn}

The compact open topology on the space $\FF$ of self-maps 
of a dotted metric space is defined by taking the family of sets
$U_K (f) = \{ g \in \FF: g(x) = f(x), ~~~ \forall x \in K; K \subset X ~~~\mbox{finite}  \}$ as a basis for the topology on $\FF$.

The following easy 
 observation will turn out to be quite useful. Note that we do not need any extra geometric assumptions (e.g. hyperbolicity)
on $X$ in the Lemma below.

\begin{lemma} Let $(X,d)$ be a proper dotted geodesic metric space.
Let $L$ be a closed subset of $\FF$, the collection  of  self-maps 
of $X$ equipped with the compact open topology. Further suppose that there exist $K\geq 1, C, \epsilon \geq 0$ and  $x \in X$ such that
for all $g \in L$, $g$ is a
 $(K, \epsilon)$ quasi-isometry of $X$ and
 $d(x, g(x)) \leq C$. Then $L$ is compact. Hence any group of
uniform quasi-isometries of a proper dotted metric space $X$ is locally compact.
\label{lc}
\end{lemma}

\noindent {\bf Proof:} Since $N_C(x)$ is finite (by properness), it suffices to prove that 
for all $y \in N_C(x)$,
$\{ g \in L \colon g(x) = y \}$ is compact. 

Hence, without loss of generality assume that there exists $y \in X$ such that $g(x) = y$ for all $g \in L$.
Since $X$ is proper, $X$ is countable. Let $X = \{ x = x_1, x_2, \cdots , x_n, \cdots \}$ be an enumeration of the elements of $X$.
Since each $g_\alpha \in L$ is a  $(K, \epsilon)$ quasi-isometry and $g(x) = y$ for all $g \in L$, 
then for each $x_n \in X$, there exists a finite set $K_n$ such that $g(x_n) \in K_n$ for all $g \in L$.

Then, given any infinite collection of $g_\alpha$'s in $L$, we can pass to a sequence $\{g_1, g_2, \cdots , g_n, \cdots \}$
such that for each $x_n \in X$, there exists $y_n \in X$ with $g_i(x_n) = y_n$ for all $i \leq n$. Let $g_\infty (x_n) = y_n$ for all $n$.

Then $g_n$ converges to $g_\infty$ in  $\FF$, the collection  of  self-maps 
of $X$ equipped with the compact open topology and we are done. $\Box$

\smallskip

Lemma \ref{lc} may be thought of as a coarsening of the fact that the stabilizer of a point in the isometry group
of a Riemannian manifold is compact.

For a non-elementary Gromov hyperbolic group $G$ we shall construct a certain pseudo-metric space which will come in handy.
It is known \cite{gromov-hypgps} that $G$ acts cocompactly on the collection $\partial^3 G$ of distinct triples on the boundary $\partial G$ of $G$. 
Let $K$ be a (closed) fundamental domain for this action. Choose a point $p$ in the interior $Int (K)$ of $K$. Define 
$\rho(g(p), h(p)) = 1$ if $g(K) \cap h(K) \neq \emptyset$. Also for $x \in g(Int(K))$ define $\rho(g(p), x) = 0$. For   $x \in 
\partial^3 G \setminus \bigcup_g  g(Int(K))$, let $g_1, \cdots , g_m$ be the collection of all
elements of $G$ such that $x \in \bigcap_i  g_i(K)$. Choose one of the elements $g_1, \cdots , g_m$, say $g_i$ and define 
$\rho(g_i(p), x) = 0$ and $\rho(g_j(p), x) = 1$ for $j \neq i$. Now define a {\it dotted path metric} on 
$\partial^3 G$ by 
\begin{center}
$\rho (x, y) = inf \{ n:$ There exists a sequence $ x=x_0, g_1(p), g_2(p), \cdots , g_n(p), x_{n+1} = y$ \\
such that $\rho (x_0, g_1(p)) =
0 = \rho (x_{n+1}, g_n(p))$ \\ 
and $\rho (g_i(p), g_{i+1}(p)) = 1$ for $i=1 \cdots n-1  \}$.
\end{center}

\begin{obs} \label{milnor-s}{\rm 
The pseudo-metric space $(\partial^3 G , \rho )$ is quasi-isometric to any Cayley graph $\Gamma$ of $G$ with respect to a finite generating
set. The proof of this fact is an easy modification of the {\v S}varc-Milnor Lemma 
(see the proof of Proposition 8.19, p. 140 of \cite{br-h}). The map $\phi : \Gamma  \rightarrow (\partial^3 G , \rho )$
given by $\phi (g) = g(p)$ gives the required quasi-isometry. }\end{obs}

\subsection{Patterns}\label{pats}

\begin{defn} Let $G$ be a hyperbolic group acting geometrically (i.e. freely, cocompactly
 and properly discontinuously by isometries) on a hyperbolic metric
space $\Hyp$. 
A {\bf symmetric pattern}
 of closed convex (or quasiconvex) sets in $\Hyp$ is a $G$-invariant countable collection 
$\JJ$ of  convex (or quasiconvex) sets such that \\
1)  The stabilizer $H$ of $J \in \JJ$ acts cocompactly on $J$. \\
2) $\JJ$ is the orbit of some (any) $J \in \JJ$ under $G$. 
\end{defn}

This definition is slightly more restrictive than Schwartz' notion
of a symmetric pattern
 of geodesics, in the sense that he takes
$\JJ$ to be a finite union of orbits of geodesics, whereas Condition (2)
above forces $\JJ$ to consist of one orbit.
All our results go through with the more general definition, where
$\JJ$ is a finite union of orbits of closed convex (or quasiconvex) sets,
but we restrict ourselves to one orbit for ease of exposition.

Suppose that
$(X_1, d_1), (X_2,d_2)$ are metric spaces. Let $\JJ_1, \JJ_2$ be collections of closed subsets of $X_1, X_2$ respectively. Then $d_i$ induces a {\em pseudo-metric} (which, by abuse of notation, we continue to refer to as $d_i$) on $\JJ_i$ for $i = 1, 2$. This is just the ordinary (not Hausdorff) distance between closed subsets of a metric space. 

In particular, 
consider two hyperbolic groups $G_1, G_2$ with quasiconvex subgroups $H_1, H_2$, Cayley graphs $\Gamma_1, \Gamma_2$. Let $\LL_j$  for $j = 1, 2$ denote the collection of translates of limit sets of $H_1, H_2$ in $\partial G_1, \partial G_2$ respectively. Individual members of  the collection $\LL_j$ will be denoted as $L^j_i$.
Let $\JJ_j$ denote the collection 
$\{ J_i^j = J(L_i^j): L_i^j \in \LL_j \}$ of joins of limit sets.
Recall that the join of a limit set
$\Lambda_i$ is the union of bi-infinite geodesics in $\Gamma_i$ with end-points in $\Lambda_i$. This is a uniformly quasiconvex set and lies at a bounded Hausdorff distance from the Cayley graph of the subgroup $H_i$ (assuming that the Cayley graph of $H_i$
is taken with respect to a finite generating set which is contained in the generating set of $G_i$).
Following Schwartz \cite{schwarz-inv}, we define:

\begin{defn} A bijective map $\phi$ from $\JJ_1 \rightarrow \JJ_2$  is said to be uniformly proper  if there exists a function $f: \natls \rightarrow \natls$ such that \\
1) $d_{G_1} (J(L_i^1) ,J(L_j^1)) \leq n \Rightarrow d_{G_2} (\phi(J(L_i^1)) ,\phi(J(L_j^1))) \leq f(n)$ \\
2)  $ d_{G_2} (\phi(J(L_i^1)) ,\phi(J(L_j^1)))\leq n \Rightarrow d_{G_1} (J(L_i^1) ,J(L_j^1)) \leq f(n)$. 

When $\JJ_i$ consists 
of all singleton subsets of $\Gamma_1, \Gamma_2$, we shall
refer to $\phi$ as a uniformly proper map from 
$\Gamma_1$ to $ \Gamma_2$.
\label{uproper}
\end{defn}

The proof of the following Theorem can be culled out of \cite{mahan-relrig}.
We give a proof for completeness.
\begin{theorem} Let $H$ be an infinite quasiconvex subgroup of a hyperbolic group $G$ such that $H$ has infinite index in $G$.
Let $\Gamma$ be a Cayley graph of $G$ with metric $d$.
Let $L$ be the limit set of $H$ and $\LL$ be the collection of translates of $L$ under $G$. There exists a finite collection
$L_1, \cdots L_n$ of elements of $\LL$ such that the following holds. \\
For any $K, \epsilon$, there exists a $C$ such that if $\phi : \Gamma \rightarrow \Gamma$ is a pattern-preserving
$(K, \epsilon )$-quasi-isometry of $\Gamma$ with $\partial \phi (L_i) = L_i$ for $i = 1 \cdots n$, then $d(\phi (1), 1) \leq C$.
\label{origin}
\end{theorem}

\begin{proof} Let $J_i$ denote $J(L_i)$ for $L_i \in \LL$. Also let $B_k(1)$ denote the $k-$ neighborhood of $1 \in \Gamma$.
In \cite{mahan-relrig} p. 1706, we show that that there exists $M \in \natls$ such that for all $k\geq M$ the collection
 $$\{ J_i : B_k (1) \cap J_i \neq \emptyset \}$$ contains a pair $J_p, J_q$ such that $L_p \cap L_q =  \emptyset$.
Further (\cite{mahan-relrig} p. 1707) 
for any $K_1$, there exists $D $, such that $\{ z \in \Gamma : d(z,J_p) \leq K_1, d(z,J_q) \leq K_1 \}$ has diameter less than $D$.

Let $L_1, \cdots L_n$ be all the elements of $\LL$ such that $B_M (1) \cap J_i \neq \emptyset$. Suppose 
$\partial \phi (L_i) = L_i$ for $i = 1 \cdots n$. Then there exists $C_0= C_0(K, \epsilon )$ such that $\phi (J_i) $ lies in a $C_0$ neighborhood
of $J_i$ for  $i = 1 \cdots n$. Also $d(\phi (1), \phi (J_i)) \leq MK + \epsilon$  for  $i = 1 \cdots n$. Hence 
 $d(\phi (1), J_i) \leq MK + \epsilon + C_0$  for  $i = 1 \cdots n$. 

  Choose
$K_1 = MK + \epsilon + C_0$. Then $d(z,J_p) \leq K_1, d(z,J_q) \leq K_1$ for $z=1$ or $z= \phi (1)$. Hence $d(1, \phi (1)) \leq D$ where $D$ is 
the real number determined by the last assertion in the first paragraph of this proof. Taking $C=D$ we are done.
\end{proof}

\begin{defn}  Let $H$ be an infinite quasiconvex subgroup of a hyperbolic group $G$ such that $H$ has infinite index in $G$.
The group $PP(G,H)$ of pattern-preserving maps for such a pair $(G,H)$  is defined as the group
of homeomorphisms of $\partial G$ that preserve the collection of translates $\LL$, i.e. $PP(G,H) = \{ \phi \in Homeo (\partial G): \phi(L) \in \LL,
~~\forall L \in \LL \}$. 
The group $PPQI(G,H)$ of pattern-preserving quasi-isometries
 for a pair $(G,H)$ as above is defined as the subgroup of $PP(G,H)$ consisting of homeomorphisms $h$
of $\partial G$ such that $h = \partial \phi$ for some quasi-isometry $\phi : \Gamma \rightarrow \Gamma$.
The topology on  $PP(G,H)$ or  $PPQI(G,H)$ is inherited (as a subspace) from the uniform topology on $Homeo (\partial G)$.
\end{defn}

\begin{prop} \cite{mahan-relrig} 
The collection $\LL$ is discrete in the Hausdorff topology on the space of closed subsets of $\partial G$, i.e. for all $L \in \LL$,
there exists $\epsilon > 0$ such that $N_\epsilon (L) \cap \LL = L$, where  $N_\epsilon (L) $ denotes an $\epsilon$ neighborhood of $L$
in the Hausdorff metric. Further, for every $\epsilon > 0$ and any visual metric $d$ on $\partial G$, the number of elements
of $\LL$ of diameter greater than $\epsilon$ is finite. \label{discrete} \end{prop}

\begin{obs} {\rm
We observe now that  $PP(G,H)$  is closed in  $Homeo (\partial G)$ equipped with the uniform topology. To see
this, assume that $f_n \in PP(G,H)$ and $f_n \rightarrow f$  in  $Homeo (\partial G)$ equipped with the uniform topology.
If $f \notin PP(G,H)$ there exists $L \in \LL$ such that $f(L) \notin   \LL$. Since $f_n$ converges to $f$ in the uniform topology,
$f_n(L) \rightarrow f(L)$ in the Hausdorff metric on closed subsets of $\partial G$. Since $f$ is a homeomorphism, there 
 exists $\epsilon > 0$ such that  the diameter of $f(L) $ is greater than $2\epsilon$
and hence  there exists  $N \in \natls$ such that  the diameter of $f(L) $ is greater than $\epsilon$
 for all $n \geq N$. Since the number of elements
of $\LL$ of diameter greater than $\epsilon$ is finite
by Proposition \ref{discrete}, then  (after passing to a subsequence if necessary) it follows
that there exists $L_1 \in \LL$  such that $f_n(L) = L_1$ for all $n \geq N$. Since $L, L_1$
are closed and $f_n \rightarrow f$  in  $Homeo (\partial G)$ equipped with the uniform topology, it follows that $f(L) = L_1 \in \LL$,
a contradiction. This shows that  $PP(G,H)$  is closed in  $Homeo (\partial G)$.

The same is true for  $PPQI(G,H)$.}
\end{obs}

{\it Henceforth, whenever we refer to  $Homeo (\partial G)$ as a topological group, we shall assume that it is  equipped with the {\bf uniform topology}.}

Combining Lemma \ref{lc} with Theorem \ref{origin}, we get 

\begin{cor} Let $H$ be an infinite quasiconvex subgroup of infinite index in a hyperbolic group $G$.
Let $\Gamma$ be a Cayley graph of $G$ with metric $d$ and $Q  \subset PPQI(G,H)$ be any group of boundary values of 
uniform quasi-isometries of (the vertex set of) $\Gamma$.
Let $L$ be the limit set of $H$ and $\LL$ be the collection of translates of $L$ under $G$. There exists a finite collection
$L_1, \cdots L_n$ of elements of $\LL$ such that $Q_0 = \cap_{i = 1 \cdots n} Stab(L_i)$ is compact, where $Stab(L_i)$ denotes the
stabilizer of $L_i$ in $Q$. (Here $Q$ is equipped with the uniform topology as a subspace of $Homeo (\partial G)$.)
\label{upp-cpt}
\end{cor}

\begin{proof} Let $K \geq 1, \epsilon \geq 0$ be such that each $q \in Q$ is the boundary value of a $(K, \epsilon)$-quasi-isometry.
Let $L_1, \cdots L_n$ be elements of $\LL$ as in Theorem \ref{origin}.
For each $q\in Q$, let $\phi_q$ be a $(K, \epsilon)$-quasi-isometry such that $\partial \phi_q =q$. Let $\FF_0 = \{ \phi_q:
\partial \phi_q \in Q_0 \}$.

By Theorem \ref{origin}, there exists $C=C(K, \epsilon)$ such that $d(\phi_q (1), 1) \leq C$ for all $\phi_q  \in \FF_0$.
Hence, by  Lemma \ref{lc}, $\FF_0$ is compact in the compact open topology on the space of self-maps of the dotted geodesic metric space consisting
of the  vertex set of $\Gamma$.

Suppose $\{ \phi_{q_i} \} \subset  \FF_0$ is a sequence converging to $\phi_q$ in the compact open topology.
Hence $q_i$ and $q$ are boundary values of $(K, \epsilon)$-quasi-isometries $\phi_{q_i}$ and $\phi_q$ respectively such that 
$\phi_{q_i}, \phi_q$ agree on larger and larger subsets as $i \rightarrow \infty$. It follows that $q_i \rightarrow q$ in
$Homeo (\partial G)$ (with the uniform topology) by stability of quasigeodesics \cite{GhH} as $G$ is hyperbolic.
\end{proof}

Also, since   $PP(G,H)$ or   $PPQI(G,H)$  is a closed subgroup  of  $Homeo (\partial G)$ equipped with
the uniform topology,  we have

\begin{cor}Let $H$ be an infinite quasiconvex subgroup of a hyperbolic group $G$ such that $H$ has infinite index.
Let $\Gamma$ be a Cayley graph of $G$ with metric $d$ and $Q \subset PPQI(G,H)$ be any group of boundary values of
pattern-preserving  uniform
quasi-isometries of (the vertex set of) $\Gamma$.
Then $Q$, equipped with the topology inherited from $Homeo (\partial G)$ is locally compact. 
\label{upp-lc}
\end{cor}

\noindent {\bf Proof:} We include a slightly different direct proof here. The collection $\LL$ is discrete by Proposition \ref{discrete}.
Consider the finite collection $L_1, \cdots L_n$ in Corollary \ref{upp-cpt}. There exists $\epsilon > 0$ such that 
$N_\epsilon (L_i) \cap \LL = L_i$ for all $ i = 1 \cdots n$. Define 

\begin{center}

$N_\epsilon (Id) = \{ q \in Q \colon d_{\bdy G} (x, q (x)) \leq \epsilon $ for all $ x \in {\bdy G} \}$

\end{center}

\noindent where $d_{\bdy G}$ denotes some visual metric on ${\bdy G}$. Then 
$N_\epsilon (Id) \subset Q_0 = \cap_{i = 1 \cdots n} Stab(L_i)$, which is compact by  Corollary \ref{upp-cpt}.
 Hence the Corollary. $\Box$

\subsection{Boundaries of Hyperbolic metric spaces and the Newman-Smith Theorem}

Let $L$ be one of the rings $\Z$ or $\Z_p$, for $p$ a prime.

\begin{defn} \label{d:gm}(\cite{bredon-sheaf}, p.329) 
 An {\em{$m$-dimensional
homology manifold over $L$}} (denoted $m$-{\rm{hm}}$_L$) is
a locally compact Hausdorff space  $X$  with finite homological dimension over  $L$, that has the local homology properties of a manifold, i.e. for
all $ x \in X$, $H_m (X, X \setminus \{ x \}) = L$ and 
$H_i (X, X \setminus \{ x \}) = 0$ for $i \neq m$.\\
Further, if $X$ is an $m$-{\rm{hm}}$_L$ and $H^c_*(X; L) \cong H^c_*( {\mathbb{S}}^m; L)$ then $X$
 is called a {\em{generalized $m$-sphere}} over $L$. (Here $H^c_*$ denotes {\v C}ech homology and $H_*$ denotes ordinary singular
homology.)
 \end{defn}
 
 For  homology manifolds, the existence of a local orientation was proven
by Bredon in \cite{bredon-mmj}.

The related notion of a 
{\em cohomology manifold over $L$}, denoted $m$-{\rm{cm}}$_L$ 
is defined by Bredon in \cite{bredon-sheaf}, p. 373.  If $L=\Z_p$,  a connected space $X$ is an $n$-{\rm{cm}}$_{L}$ if and only it if it is an $n$-{\rm{hm}}$_{L}$ and is  locally connected (\cite{bredon-sheaf}, p. 375 Theorem 16.8, footnote).  

We shall  be using the following Theorem 
which is a result that follows from work of
Bestvina-Mess \cite{bestvina-mess} and Bestvina
\cite{bestvina-mmj} (See also Swenson \cite{swenson-mmj}  , Bowditch \cite{bowditch-cutpts} and Swarup
\cite{swarup-cutpts}).

\begin{theorem} \label{bestvina00} Boundaries $\partial G$
of PD(n) hyperbolic
groups $G$ are locally connected
homological manifolds (over the integers) with the homology
of a sphere of dimension $(n-1)$. Further, if $G$
acts freely, properly discontinuously, cocompactly on a contractible complex $X$
then  $H_n^{LF}(X) = H_{n-1} \partial G$, where $H_n^{LF}$ denotes locally finite homology. 
\end{theorem}

In fact, one of the main results of
\cite{bestvina-mess} asserts that the (reduced) \v{C}ech cohomology groups of $\partial G$
vanish except in  dimension $(n-1)$.  Bestvina \cite{bestvina-mmj} also shows that the 
(reduced) Steenrod homology groups of $\partial G$ vanish except in  dimension $(n-1)$. Since 
$\partial G$ is compact metrizable, Steenrod homology coincides with \v{C}ech homology (see for instance \cite{milnor-steenrod}).
Further, for locally connected metrizable compacta such as $\partial G$, the \v{C}ech (co)-homology groups coincide
with singular (co)-homology groups (see pg. 107 of \cite{lefschetz_bk2}). Hence 
the singular (as well as \v{C}ech) homology and cohomology of $\partial G$ coincide with that of a sphere of dimension $(n-1)$,
i.e. $H_0(\partial G) = H_{n-1}(\partial G) = H^0(\partial G) = H^{n-1}(\partial G) = \mathbb{Z}$ and all other homology and cohomology
groups with $\mathbb{Z}$ coefficients are zero.
Finally, using the Universal Coefficient Theorem for homology (Theorem 3A.3 of \cite{hatcher-bk})
$H_i(\partial G ; {\mathbb{Z}}_p)   = H_i(\partial G) \otimes
 {\mathbb{Z}}_p$ since all the integral homology groups are torsion-free. Similarly for cohomology groups.
Thus we have the following strengthening of Theorem \ref{bestvina00}:

\begin{theorem} \label{bestvina} Let $L$ denote $\Z$,
the integers  or $\Z_p$, the integers mod $p$. Boundaries $\partial G$
of PD(n) hyperbolic
groups $G$ are locally connected
(co)homological manifolds (over $L$) with the (co)homology
of a sphere of dimension $(n-1)$. 
\end{theorem}

We shall be using Theorem \ref{bestvina} in conjunction with the following Theorem of Newman and Smith, which as stated below
is a consequence of the work of several people (see below).

\begin{theorem} (Newman \cite{newman}, Smith \cite{smith3})
Let $(X,d)$ be a compact $\Z_p$-cohomology manifold for all $p$, having finite topological (covering) dimension, and 
equipped with a metric $d$. There exists $\epsilon > 0$ such that (for any $n$) if $\Z_n$
acts effectively on $X$, then the diameter of some orbit is greater than
$\epsilon$.
\label{ns}
\end{theorem}

Newman proved the above Theorem for closed orientable manifolds \cite{newman}.
Smith \cite{smith3} generalized it to locally compact spaces satisfying certain homological
regularity properties.
Building on work of Yang \cite{yang},
Conner and Floyd (\cite{conner-floyd}, Corollary 6.2)
proved that cohomological manifolds of finite
topological (covering) dimension satisfy the regularity properties required by Smith's 
theorem. (The theorem also holds for a somewhat more general class of spaces, called
`finitistic spaces' by Bredon in \cite{bredon-tg}, but we shall not require this).
For a historical account of the development of the theory of generalized/homology manifolds and their connection with Smith manifolds see \cite{raymond}.

We shall also be using a theorem on homological consequences
of actions of $p$-adic transformation groups on homology manifolds.

\begin{theorem} (Yang \cite{yang-p})
Let $X$ be a compact homology $n$-manifold admitting an effective  $K$-action, where $K=Z_{(p)}$  is the group
of $p$-adic integers. Then the homological dimension of $X/K$ is $n+2$.
\label{yang}
\end{theorem}

\section{Ahlfors Regular Metric Measure Spaces}

\subsection{General Facts} 
We refer to  \cite{heinonen-bk} and \cite{bonk-kleiner}  for details on quasi-symmetric maps and
metric measure spaces.

\medskip

\noindent {\bf Quasi-symmetric Maps and Quasi-isometries}

\begin{defn}
    Let $ f : X \rightarrow Y$ be a homeomorphism between metric spaces $(X, d_X )$ and $(Y, d_Y )$.

Then $f$ is  {\bf $\eta$-quasi-symmetric } for some homeomorphism $\eta : [0, \infty ) \rightarrow [0, \infty )$ if
$$\frac{d_Y(f(x_1),f(x_2))}{d_Y(f(x_1),f(x_3))}
\leq \eta\left (\frac{d_X(x_1,x_2)}{d_X(x_1,x_3)}\right )$$
for every triple $(x_1,x_2,x_3)$ of distinct points in $X$. \end{defn}

The inverse of a quasi-symmetric map is 
also quasi-symmetric. In fact the right generalization of quasiconformal maps in $\R^n$ to metric measure spaces
are quasi-symmetric maps (cf. Ch. 10, 11 of \cite{heinonen-bk}).

\begin{prop} (\cite{heinonen-bk} p. 79) If $f: X \rightarrow Y$ is an
$\eta$-quasisymmetric homeomorphism, then $f^{-1}: Y \rightarrow X$ is an
$\eta_1$-quasisymmetric homeomorphism, where $\eta_1(t) = 1/\eta^{-1}(t^{-1})$ for $t > 0$. Further, if 
$f: X \rightarrow Y$ and $g: Y \rightarrow Z$ are
$\eta_f$ and $\eta_g$-quasisymmetric homeomorphisms respectively, then $g \circ f: X \rightarrow Z$
is $\eta_g \circ \eta_f$-quasisymmetric. \label{qs-gp} \end{prop}

\begin{rmk} {\rm It follows from the definition of $\eta$-quasi-symmetry that $d_X(x_1,x_2) \leq d_X(x_1,x_3) \Rightarrow 
d_Y(f(x_1),f(x_2)) \leq \eta(1) d_Y(f(x_1),f(x_3))$.  Now let $B= B(x,r) \subset X$ denote the closed ball
of radius $r$ about $x \in X$ and $\partial B = \{ u \in B: d_X(x,u) = r \}$.  Hence
if $X$ is connected and $y, z \in \partial B $, then $d_Y(f(x),f(y)) \leq \eta(1) d_Y(f(x),f(z))$.} \label{ball} \end{rmk}

\begin{lemma} Let $g: X \rightarrow X$ be an $\eta$-quasi-symmetric map of a compact
metric space $(X,d)$ to itself. 
Then $g^{-1}$ is an $\eta_1$-quasi-symmetric map by Proposition \ref{qs-gp}. Let $c=\eta_1^{-1} (1)$.
Let $B = B(x,r) = \{ y \in X: d(x,y) \leq r \}$ be the closed ball of radius $r$ about $x$.
Let $\partial B = \{ y \in X: d(x,y) = r \}$ denote the sphere of radius $r$ about $x$. Assume that $\partial B \neq \emptyset$. Let $s = d(g(x), g(\partial B))$.
Then $B(g(x), cs) \subset g(B(x,r)$. \label{nag} \end{lemma}

\noindent {\bf Proof:} It follows from compactness of $X$
that $\partial B$ is compact. Therefore $g(\partial B)$ is compact (and non-empty by hypothesis). 
Hence there exists
 $y_0 \in \partial B$  such that $d(g(x), g(y_0)) = s$. 
Suppose $\frac{d(g(x),w)}{d(g(x),g(y_0))} \leq c$ for some $w \in X$. Then 
$\frac{d(x,g^{-1}(w))}{d(x,y_0)} \leq \eta_1(c) = 1$. Hence $d(x,g^{-1}(w)) \leq d(x,y_0) = r$.
Hence $w \in g(B)$. $\Box$

\smallskip

It is a standard fact that the boundary values of  quasi-isometries of proper hyperbolic metric spaces are quasi-symmetric
maps for any visual metric on the boundary:

\begin{lemma} \label{qiqs} (\cite{BS} Theorem 5.2.17, p. 55; \cite{vaisala-notes} Theorem 5.35)\\
Suppose that $X$ and $Y$ are $\delta$-hyperbolic spaces equipped with base-points and that
$f : X \rightarrow Y$ is a base point preserving $(\lambda , \mu)$-quasi-isometry. Let $(\partial X, \rho_X)$ 
and $(\partial Y, \rho_Y)$ be their respective boundaries equipped with visual metrics.   Then $f$ extends to an 
$\eta$-quasisymmetric homeomorphism $\partial f :(\partial X, \rho_X) \rightarrow (\partial Y, \rho_Y)$
     with $\eta$
depending only on $\delta ,\lambda , \mu$. \end{lemma}

A converse result follows from work of Paulin \cite{paulin-jlms}. We state it in the form we shall need it. 

\begin{theorem}   Suppose that $X$ is (the Cayley graph with respect to a finite generating set of) 
 a non-elementary hyperbolic group.
Further suppose that $X$ is equipped with a base-point. Let $Q$
be a compact group of $\eta$-quasisymmetric homeomorphisms in the uniform topology on $\partial X$. Then there exist 
$(\lambda , \mu)$ such that each $q \in Q$ may be realized as the boundary value of a $(\lambda , \mu)$-quasi-isometry
of $X$ fixing the base-point. \label{qiqs2} \end{theorem}

\noindent {\bf Proof:} By Paulin's work \cite{paulin-jlms}, there exist $\lambda_1 , \mu_1$ such that every element 
$q \in Q$ may be realized as the boundary value of a $(\lambda_1 , \mu_1)$-quasi-isometry of $X$. 
 We denote the
map induced by $q$ on $(\partial^3 X, \rho)$  by $q^3$. Here $(\partial^3 X, \rho)$ is the pseudo-metric space in
Observation \ref{milnor-s}.

Since $X$ is quasi-isometric
to $(\partial^3 X, \rho)$ for the dotted path-metric constructed on the set of distinct triples (Observation \ref{milnor-s}), 
 every element 
$q \in Q$ may be realized as the boundary value of a uniform quasi-isometry $q^3$ of $(\partial^3 X, \rho)$.

Since $Q$ is compact, any $q^3 \in Q$  maps a fixed triple $(a_1, a_2, a_3) \in (\partial^3 X, \rho)$ to triples that are uniformly
separated in $(\partial X, d)$ where $d$ denotes a visual metric, i.e. there exists $\epsilon > 0$ such that
for all $q \in Q$,  $d(q(a_i), q(a_j)) > \epsilon$ for $i \neq j$. Hence $ \rho ((a_1, a_2, a_3), q^3(a_1, a_2, a_3))$ is
uniformly bounded. Choose $(a_1, a_2, a_3) \in (\partial^3 X, \rho)$ as the base-point in 
$(\partial^3 X, \rho)$.

 Define a map
$\phi_q :  (\partial^3 X, \rho) \rightarrow  (\partial^3 X, \rho)$ as follows:\\
1) $\phi_q$ maps $q^3(a_1, a_2, a_3)$ and all points at $\rho$ distance zero from it to $(a_1, a_2, a_3)$.\\
2)  $\phi_q$ fixes all other points of  $(\partial^3 X, \rho)$. \\

Then $\phi_q \circ q^3$ is a uniform (independent of $q$) quasi-isometry of $(\partial^3 X, \rho)$ fixing a base-point. 
Post-composing further by the inverse of the quasi-isometry from $X$ to $(\partial^3 X, \rho)$ and moving back a base-point in $X$ again
if necessary, we get the required result. $\Box$

\medskip

\noindent {\bf Ahlfors Regular and Doubling Spaces} 

\begin{defn} \label{doubling-def}
A metric space $X$ is said to be {\it doubling} if  for all $\lambda \geq 1$ there exists $N \in \mathbb{N}$ such that
for all $x \in X$ the ball $\lambda B(x,r) = B(x, \lambda r)$ can be covered by $N$ balls of radius $r$. \end{defn}

We shall need
a special case of a Theorem of Bonk and Schramm \cite{bonk-schramm}.

\begin{theorem} {(\bf Bonk-Schramm \cite{bonk-schramm})}
Let $X$ be a Gromov-hyperbolic group or a complete simply connected Hadamard manifold of pinched
negative curvature. Then $(\partial X, d)$ is doubling for any visual
metric $d$.
\label{doubling}
\end{theorem}

\begin{defn} \label{haus-def}  Let $(X,d)$ be a  metric  space. Then the $k-$dimensional Hausdorff measure of $X$ is defined by
\begin{center}
$inf \{  \sum_i r_i^k: $ there exists a  cover of $X$ by balls $B_i$  of radius $r_i \}$
\end{center}
The 
Hausdorff dimension  $dim_{haus}(X)$ of $X$ is defined to be the infimum of all $k \geq 0$ such that the $k-$dimensional Hausdorff measure of $X$
is zero.
\end{defn}

In this paper we shall have occasion to use the   easier to compute notion of  Minkowski dimension where the balls used to cover $X$
are of equal radii. See \cite{sullivan-pihes, coornaert, bj} for the equivalence in
 our situation (also see \cite{falc} for  very general sufficient conditions ensuring
equality of Hausdorff and Minkowski dimensions).  

\begin{defn} \label{ahlforsreg-def}
Let $(X,d, \mu )$ be a compact metric measure space, i.e. a metric space equipped with a Borel measure. We say that
$X$ is Ahlfors $Q$-regular, if there exists $C_0 \geq 1$ such that for all $0\leq r\leq dia(X)$, and any ball
$B_r(x) \subset X$, the measure $\mu (B_r(x))$ satisfies $\frac{1}{C_0} r^Q \leq \mu (B_r(x)) \leq C_0r^Q$. \end{defn}

When
$Q$ is omitted in Definition \ref{ahlforsreg-def}, we assume that $Q$ is the Hausdorff dimension and $\mu$ the Hausdorff measure.

The relevance
to the present paper comes from the following.

\begin{theorem} {(\bf Coornaert \cite{coornaert})}
Let $G$ be a hyperbolic group. Then $(\partial G, d, \mu)$ is Ahlfors regular for any visual
metric $d$ and the associated Hausdorff measure $\mu$.
\label{ahlfors}
\end{theorem}

Certain very general conditions ensure  Ahlfors regularity. 
 Corollary 14.15 of \cite{heinonen-bk} asserts that a metric space is quasi-symmetrically equivalent
to an Ahlfors regular space if and only if it
is {\it uniformly perfect} and carries a {\it doubling measure}. We shall not be  needing the precise definitions
of  these terms. Suffice to say that a metric space carries a doubling measure if and only if
it is  doubling \cite{lu-sak} \cite{coif-weiss}; and connected sets are uniformly
perfect.  However for our purposes, the {\it proof} of the more restrictive Theorem 5.4 of \cite{coornaert} for hyperbolic
groups which in turn is modelled on Sullivan's work \cite{sullivan-pihes} suffices  to
ensure the following. (See Theorem 25 of \cite{sullivan-pihes} 
in particular which devotes special attention to 
non-compact finite volume manifolds or Theorem 0.2 of \cite{paulin-dg}.)

\begin{theorem}
Let $X$ be the universal cover of a complete finite volume manifold of pinched
negative curvature. Then $(\partial X, d, \mu)$ is Ahlfors regular for any visual
metric $d$ and the associated Hausdorff measure $\mu$.
\label{ahlfors2}
\end{theorem}

The following notion is essentially due to Bourdon and Pajot cf. \cite{bp}.
\begin{defn} \label{acd} Let $\Gamma$ be a hyperbolic group or more generally a hyperbolic metric space.
 We define the Ahlfors regular conformal dimension of the boundary of $\Gamma$ as follows.
\begin{center}
$ACD(\partial \Gamma ) = inf \{Q : \partial \Gamma$ admits a visual metric with Hausdorff dimension $Q \}$.
\end{center}
\end{defn}

In this paper we shall repeatedly make the assumption that the boundary $\partial \Gamma$ of 
$\Gamma$ admits a visual metric $d$
such that the Hausdorff dimension $dim_{haus}$ of $(\partial \Gamma, d)$ is less than $(dim_t + 2)$, where $dim_t$ denotes 
topological dimension. Note that by Definition \ref{acd} above, this is equivalent to saying  that $ACD (\partial \Gamma) < 
(dim_t (\partial \Gamma) + 2) $. 

Also note that all visual metrics are quasi-symmetrically conjugate to each other. The identity
map on the underlying topological space gives the required quasi-symmetry. More generally,
for a metric space $X$, $ACD(X)$ is usually defined \cite{bp} by
\begin{center}
$ACD(X) = inf \{Q :$ there exists an Ahlfors $Q$-regular metric space $ Y$ that is quasisymmetrically equivalent to $X \}.$
\end{center}
We shall not have need for this more general notion in this paper.

\smallskip

We collect together certain notions and facts from \cite{bonk-kleiner} and \cite{heinonen-bk} 
that we shall have need for in what follows.

\begin{lemma} (\cite{heinonen-bk} Theorem 1.2, p. 2) Every family $\mathcal F$ of balls of uniformly bounded diameter in a metric space
$X$ contains a pairwise disjoint subfamily $\mathcal G$ such that $\bigcup_{B \in \mathcal F} B \subset \bigcup_{B \in \mathcal{G}} 5B$, where
$\lambda B$ denotes a ball concentric with $B$ and radius $\lambda$ times that of $B$. \label{covering} \end{lemma}

\begin{cor} Let $X$ be a doubling metric space. There exists an $M$ such that the following holds. \\ Let  ${\mathcal F}_r$ be the family of all balls
of radius $r$ in $X$. Let ${\mathcal G}_r$ be any  pairwise disjoint subfamily such that $X = \bigcup_{B \in {\mathcal{G}}_r} 5B$.
Then for any $x \in X$, the cardinality of the set $\{ B \in {\mathcal{G}}_r : x \in 5B \}$ is less than or equal to $M$. \label{covercor} \end{cor}

\noindent {\bf Proof:} We omit the suffix $r$ for convenience.
The existence of  a  pairwise disjoint subfamily  $\mathcal G$ of  $\mathcal F$
such that $X = \bigcup_{B \in \mathcal{G}} 5B$ is guaranteed by Lemma \ref{covering}. \\
We now argue by contradiction. Suppose that no such $M$ exists. Then for all $N \in \mathbb{N}$, there exists $r_N > 0$, 
a pairwise disjoint subfamily $\mathcal G$ of balls of radius $r_N$ in $X$, and an $x \in X$ such that 
the cardinality of the set $\{ B \in \mathcal{G} : x \in 5B \}$ is greater than $N$. Hence the $6r_N$ ball around $x$ contains
more than $N$ distinct (and hence disjoint) elements of $\mathcal G$. Since $N$ is arbitrary, this contradicts the hypothesis that $X$ is doubling. $\Box$

\medskip

We start with a  general result about compact group actions on 
doubling, Ahlfors regular compact 
metric measure spaces, e.g. boundaries of one-ended hyperbolic groups.
\begin{lemma}
Let $(X , d, \mu )$ be a connected, doubling, Ahlfors regular compact 
metric measure space having Hausdorff dimension $\QQ \geq 1$ (for instance the boundary of a one-ended hyperbolic group
or the boundary of the universal cover of a finite volume manifold of pinched negative curvature by Theorems \ref{doubling}
\ref{ahlfors} and \ref{ahlfors2}). Let $K$ be a compact
topological group acting by   uniformly $\eta$-quasi-symmetric 
 maps on $(X, d, \mu )$ and equipped with a Haar
measure of unit mass. Let $C = \eta (1)$, $c=\eta_1^{-1} (1)$, and let $d_K$ be the average
metric on $X$ given by $d_K(x,y) = \int_K d(g(x),g(y)) dg$. Then the Hausdorff dimension
of $(X, d_K, \mu )$ does not exceed $\QQ$. \label{hausdim} \end{lemma}

\noindent {\bf Proof:}   We adapt an argument of Repovs-Scepin \cite{repovs} and Martin \cite{martin-hs}
to the present context. Assume after normalization, that the total Hausdorff measure of $X$ is one.

Cover $X$ by a family  of balls of radius $r>0$ measured with respect to $d$. By Lemma
\ref{covering} we can choose a family $\BB_r$ of pairwise disjoint balls 
of radius $r$ such that $X = \bigcup_{B \in \BB_r} 5B$. 
As in \cite{martin-hs}, it is enough
 to find a uniform
bound (independent of $r$) 
for the sum $\Sigma_{B \in \BB_r} (dia_{d_K} 5B)^\QQ$
where $dia_{d_K} 5B$ 
is the diameter of $5B$ in the invariant metric $d_K$. 

 For each $B \in \BB_r$ let  $c_B$ denote its center and let $z_B$ be a point
on the boundary of $5B$ with $dia_{d_K} B \leq 2d_K (c_B , z_B )$. (Such a point $z_B$ exists as $X$ is connected.) Then \\

\smallskip

\noindent $\Sigma_{B \in \BB_r} (dia_{d_K} 5B)^\QQ $\\

\smallskip

\noindent$\leq  2^\QQ \Sigma_{B \in \BB_r} (\int_K d(g(c_B),g(z_B)) dg )^\QQ$ \\

\smallskip

\noindent $\leq  2^\QQ  \Sigma_{B \in \BB_r}\int_K d(g(c_B),g(z_B))^\QQ dg $ (by Holder's inequality since $\QQ\geq 1$ and the
normalization condition that the total Haar measure of $K$ is one)\\  

\smallskip

\noindent                   
$\leq  (2C)^\QQ  \Sigma_{B \in \BB_r}\int_K inf_{{y_B \in \partial (5B)}} d(g(c_B),g(y_B))^\QQ dg $  
(since $K$ acts by uniformly $\eta$-quasi-symmetric
maps, and since $\partial (5B) \neq \emptyset$ by connectedness of $X$, and by Remark \ref{ball}) \\

\smallskip

\noindent
$= (\frac{2C}{c} )^\QQ \Sigma_{B \in \BB_r}\int_K (c~~ inf_{{y_B \in \partial B}} d(g(c_B),g(y_B)))^\QQ dg $.\\

\smallskip

\noindent
But by Lemma \ref{nag}, the ball of radius $c~~ inf_{{y_B \in \partial B}} d(g(c_B),g(y_B))$
is contained in $g(5B)$. Hence\\

\smallskip

\noindent
$(\frac{2C}{c} )^\QQ \Sigma_{B \in \BB_r}\int_K (c~~ inf_{{y_B \in \partial B}} d(g(c_B),g(y_B)))^\QQ dg $\\
 
\smallskip

\noindent $\leq C_0 (\frac{2C}{c} )^\QQ   \Sigma_{B \in \BB_r} \int_K \mu (g(5B)) dg $  (where $C_0$ is the constant appearing in the
definition of   Ahlfors $\QQ$-regularity of $X$)\\ 

\smallskip

\noindent $= C_0(\frac{2C}{c} )^\QQ  \int_K \Sigma_{B \in \BB_r} \mu (g(5B)) dg $  \\

\smallskip

\noindent
$\leq  C_0M(\frac{2C}{c} )^\QQ   \int_K \mu (g(X)) dg $ (where $M$ is the constant appearing in Corollary \ref{covercor} giving an upper
bound for the local multiplicity of the cover $\{ 5B: B \in \BB_r \}$) \\

\smallskip

\noindent $= C_0M(\frac{2C}{c} )^\QQ   \int_K \mu (X) dg$ (since $g$ is a homeomorphism from $X$ onto itself) \\

\smallskip

\noindent $= C_0M(\frac{2C}{c} )^\QQ$ (since $K$ is equipped with a Haar measure of unit mass and by the
normalization condition that the total Hausdorff measure of $X$ is one)\\

\smallskip

\noindent This establishes  a uniform
bound (independent of $r$) 
for the sum $\Sigma_{B \in \BB_r} (dia_{d_K} 5B)^\QQ$ and completes the proof. $\Box$

\begin{theorem} \label{hr1} (See 
\cite{hewitt-ross1} ,  Ch II, Thms 7.1 - 7.6 p. 60-61). \\
Let $G$ be a  topological group. Then there exists an exact sequence $$1 \rightarrow G_0 \rightarrow G \rightarrow H \rightarrow 1$$
with $G_0$ the connected component of the identity in $G$ and $H$  totally disconnected. If $G$, and hence
$H$, is locally compact, then $H$ contains arbitrarily small compact open subgroups (i.e. for every neighborhood
$U$ of the identity in $H$, there is a compact open subgroup $K$ contained in $U $).
\end{theorem}

Moreover the structure of $G_0$ is well-known thanks to Montgomery and Zippin (\cite{mz}, Thm 4.6) 
and Gleason \cite{gleason}.

We state a definition first.

\begin{defn} A topological group is said to have {\bf arbitrarily small torsion elements} if  for every
neighborhood $U$ of the identity, there exists an element $g \in U$, $g \neq 1$, and a positive integer $n$ such that
$g^{n} = 1$, and furthermore $g^m \in U$ for all $m \in \mathbb{N}$. 

A topological group is said to have {\bf small subgroups} if for every
neighborhood $U$ of the identity, there exists a  subgroup $H \neq \{ 1 \}$ such that $H$ is contained in $U$.
\end{defn}

The following is a celebrated Theorem of Montgomery-Zippin and Gleason.

\begin{theorem} (Sections 4.5, 4.6, 4.9 of  \cite{mz}, \cite{gleason}) 
Let $G_0$ be    the connected component of the identity in a locally compact topological group $G$
such that $G/G_0$ is compact.
For each neighborhood $U$ of the identity in $G_0$, there exists a compact normal  subgroup $K \subset U$ such that the quotient group
$G_0/K$ is a Lie group.  \\ 
Small subgroups of connected  locally compact finite dimensional
groups are totally disconnected and belong to the center.  \\
Any locally compact finite dimensional
group with no small subgroups is a Lie group.
\label{smallsubgp}
\end{theorem}

\begin{rmk} \label{smallsubgp1} {\rm The hypothesis in Theorem \ref{smallsubgp} that $G/G_0$ is compact is superfluous. Let $H=G/G_0$
and $q : G \rightarrow H$ be the quotient homomorphism.
By Theorem \ref{hr1}, $H $ contains a compact open totally disconnected subgroup $H_1$. Let $G_1 = q^{-1} (H_1)$. Then $G_0$ is the
connected component of the identity of $G_1$ and $G_1/G_0$ is compact. Theorem \ref{smallsubgp} now applies.}\end{rmk}

We also have the following Theorem of Yamabe \cite{yamabe}. 

\begin{theorem} {\bf (Yamabe \cite{yamabe})}
A group with a neighborhood $U$ of the identity which does not
contain any non-trivial normal subgroup has a neighborhood $V$ of the identity which does
  not contain  any non-trivial subgroup.
\label{smallsubgp2}
\end{theorem}

\begin{defn} A topological group containing a cyclic dense subgroup is said to be {\bf monothetic}.
\end{defn}

It is easy to see that a monothetic group is abelian (using continuity of multiplication).
The following is a  consequence of a
structure theorem for $0$-dimensional compact monothetic groups (See \cite{hewitt-ross1}, Thm. 25.16, p. 408.)

\begin{theorem}
Any infinite $0$-dimensional (i.e. totally disconnected), compact, monothetic group $K$ contains
a copy of the $\bf{a}$-adic integers $A_{\bf a}$, where ${\bf a} = \{ a_1, a_2, \cdots \}$ is a sequence of integers
${a_i} > 1$. Hence $K$ must contain arbitrarily small torsion elements or a copy
of the group $Z_{(p)}$ of $p$-adic integers.
\label{div}
\end{theorem}

\subsection{The Hilbert-Smith Conjecture}
We are now in a position to prove the following.

\begin{theorem} \label{hs} 
Let 
 $(X , d, \mu )$ be a connected, doubling, Ahlfors regular compact 
metric measure space. Further suppose that 
$X$  is a $\Z_p$-cohomology manifold for all $p$; 
and $1 \leq dim_{haus}(X) < dim_h(X) +2$, where $dim_{haus}$ is the
Hausdorff dimension and $dim_h$ is the
homological dimension. Then $(X, d, \mu )$ does not admit an effective $Z_{(p)}$-action  by
{\bf uniform} quasi-symmetric maps, where 
$Z_{(p)}$ denotes the $p$-adic integers. 
Hence any finite dimensional locally compact group acting effectively
on $(X, d, \mu )$ by uniform quasi-symmetric maps is a Lie group.
\end{theorem}

\noindent {\bf Proof:} Let $K = Z_{(p)}$ be the compact group of $p$-adic integers
acting effectively on $X$ by uniform quasi-symmetric maps. Let $d_K$ be the average
metric on $X$ given by $d_K(x,y) = \int_K d(g(x),g(y)) dg$. By Lemma
\ref{hausdim} the Hausdorff dimension
of $(X, d_K, \mu )$ does not exceed $dim_{haus}(X)$. Then $K$ acts on $(X, d_K)$ by isometries. Hence the
 the orbit space $X/K$ admits the well-defined metric $\rho([x],[y]) = d_K (K(x), K(y))$, where 
$[x], [y]$ denote the images of $x, y$ under the quotient map by $K$.
Let $P: X \rightarrow  X/K$ be the natural quotient map. Since $P$ is clearly $1$-Lipschitz, it cannot
increase Hausdorff dimension. Hence the Hausdorff dimension of  $X/K$ is at most equal to $dim_{haus}(X)$, the 
Hausdorff dimension of  $X$, which in turn is less than $dim_h(X) +2$. Since topological dimension is majorized by
Hausdorff dimension and homological dimension is majorized by topological dimension, it follows that
the homological dimension of $X/K$ is  less than $dim_h(X) +2$.  This directly contradicts Yang's Theorem
\ref{yang} which asserts that the homological dimension of $X/K$ is  equal to $dim_h(X) +2$
and establishes the first part of the theorem.

The last statement follows from the first by standard arguments (see \cite{repovs} or \cite{martin-hs} for instance).
We outline the argument for the sake of completeness. 

Let $G$ be any finite dimensional
locally compact group acting effectively on $X$
by uniform quasi-symmetric maps. By the last statement of  Theorem \ref{smallsubgp}, if $G$ has no small subgroups,
then $G$ must be a Lie group. 

We now proceed by contradiction. If possible, let $G$ 
have small subgroups. By local compactness, we may assume that $G$ has compact small subgroups.
Hence by Yamabe's Theorem \ref{smallsubgp2},
$G$ has a sequence $(K_i)_i$ of  compact normal small subgroups such that $K_i \neq \{ 1 \}$ for all
$i$ and such that $\bigcap_i K_i = \{ 1 \}$. Let $L$ be the connected component of the identity of
$G$. Then either $K_i \cap L \neq \{ 1 \}$  or $K_i$ is totally disconnected. If  $K_i \cap L \neq \{ 1 \}$  then 
 $K_i \cap L$ is a compact normal small subgroup of $L$. Hence again by the second statement of  Theorem \ref{smallsubgp},
$K_i \cap L$  is totally disconnected. In any case, $G$ has 
a sequence $(K_i)_i$ of  non-trivial compact normal totally disconnected small subgroups.

  If $K$ is infinite,
it must contain a copy of the
$p$-adics (Theorem \ref{div}) or have arbitrarily small torsion elements. By Theorem \ref{ns},
$K$ cannot have arbitrarily small torsion elements. Hence 
$Z_{(p)}$ 
acts effectively on $X$ by uniform quasi-symmetric maps, contradicting the first assertion of the Theorem 
proved above.
$\Box$

\smallskip

Since quasi-isometries of a hyperbolic group $G$ act by quasi-symmetric maps on the boundary
$(\partial G, d)$, where $d$ is a visual metric, we have the following by 
combining Theorems \ref{bestvina}, \ref{ahlfors} with Theorem \ref{hs}.

\begin{cor} \label{hscor}
Let $\Gamma$ be a Poincare duality hyperbolic group and $Q$ be a group of (boundary values of)
quasi-isometries of $\Gamma$. $Q$ is equipped with the uniform topology. Suppose $d$ is a visual metric on $\partial \Gamma$ with 
$ dim_{haus} < dim_t +2$, where $dim_{haus}$ is the
Hausdorff dimension and $dim_t$ is the
topological dimension of $X= (\partial \Gamma,d)$. 
Equivalently suppose that $ACD(X) < dim_t + 2$ where $ACD(X)$ denotes the Ahlfors regular conformal dimension of $X$.
Then $Q$ cannot contain a copy of 
 $Z_{(p)}$, where 
$Z_{(p)}$ denotes the $p$-adic integers. Hence if $Q$
 is finite dimensional locally compact, it must be a Lie group.
\end{cor}

\noindent {\bf Proof:} By Theorem \ref{bestvina}, 
$X$  is a $\Z_p$-cohomology manifold for all $p$.   By Theorems \ref{doubling},  \ref{ahlfors} and by one-endedness of $\Gamma$,
$X = (\partial \Gamma, d, \mu )$ is a connected, doubling, Ahlfors regular compact 
metric measure space. 

The assumption of Theorem \ref{hs} that
 $1 \leq dim_{haus}$ is superfluous here as the only PD group that violates
this hypothesis is virtually cyclic, when the Corollary is trivially true.

  If possible, let $K = Z_{(p)}$ be the compact group of $p$-adic integers
acting effectively on $X$ by quasi-symmetric maps. All that remains to be shown is that 
we can extract an action of $Z_{(p)}$ by {\bf uniform}  quasi-symmetric maps.

Let $(\partial^3 X, \rho )$ be the pseudo-metric space in Observation \ref{milnor-s}.
Fix a base-point $x \in \partial^3 X$. Since $K$ is compact, $\rho (k^3 (x), x)$ is uniformly bounded.

For each $k\in K$ there exist $\lambda_k \in \mathbb{N}$ 
such that $k^3: (\partial^3 X, \rho ) \rightarrow (\partial^3 X, \rho )$ is a $(\lambda_k, \lambda_k)$-quasi-isometry. 
Note further that $q_n \rightarrow q$ in the uniform topology implies that $q_n^3 \rightarrow q^3$ in the compact open topology.

Let $U_i = \{ k \in K, \lambda_k \leq i\}$. Then $K = \bigcup_i \overline{U_i}$ is the union
of a countable family of closed sets $\overline{U_i}$. Note further that the limit of a sequence of
$(\lambda_k, \lambda_k)$-quasi-isometries is (at most) a $(\lambda_k + 2, \lambda_k + 2)$ quasi-isometry. 
By the Baire category
theorem there is some $U_C$ with nonempty interior. Translating by an element $h$ of
$K$ in the interior of $U_C$,  we may assume that $U_c$ (for some $c$ depending on $C$ and the quasi-isometry constant of $h$)
contains  
the identity. But any neighborhood of the identity in $K$ contains an isomorphic copy of $K$.

Hence we have an action of $K$ (replacing the original group by the isomorphic copy contained in the
above neighborhood of the identity) on $(\partial^3 X, \rho )$ by $(c,c)$-quasi-isometries. Since $(X = \partial \Gamma, d)$ is a visual
boundary for $(\partial^3 X, \rho )$ and since each element in $K$ moves the base-point $x$ by a uniformly
amount,  we have an action of $K$ on $(X, d, \mu )$
by {\it uniformly} quasi-symmetric maps by Lemma \ref{qiqs}. Theorem \ref{hs} now furnishes the required conclusion. $\Box$

\smallskip

\noindent {\bf Note:} 
The hypothesis $dim_{haus} < dim_t +2$ (or $ACD(X) < dim_t +2$) is clearly true for (uniform lattices in) real hyperbolic space, where
$dim_{haus} = dim_t$ as well as complex hyperbolic space, where
$dim_{haus} = dim_t+1$. Amongst rank one symmetric spaces these are the only ones of interest in the context of
pattern rigidity as quaternionic hyperbolic space and the Cayley plane are quasi-isometrically rigid
in light of Pansu's fundamental result \cite{pansu-rig}. 

\section{Pattern Preserving Groups as Topological Groups}

\subsection{Infinite Divisibility }\label{infdivs}
We begin with some easy classical facts about topological groups.  A property we shall be investigating
in some detail is the notion of `topological infinite divisibility'. The notion we introduce is weaker
than related existing notions in the literature.

\begin{defn} An element $g$ in a topological group $G$ will be called
 {\bf  topologically infinitely divisible}, if there exists a sequence
of  symmetric neighborhoods $(U_k)_k$ of the identity, such that $\bigcap_k U_k = \{ 1\}$ and
$g \in \bigcup_{n=1}^\infty U_k^n \subset G$ for all $k$.
 Similarly, a subgroup $H$ of $G$ is said to be {\bf  topologically infinitely divisible}, 
if there exists a sequence
of  symmetric neighborhoods $(U_k)_k$ of the identity, such that $\bigcap_k U_k = \{ 1\}$ and
$H \subset \bigcup_{1}^\infty U_k^n \subset G$ for all $k$.\end{defn}

Let $\Gamma$ be a hyperbolic group and $H$ an  infinite quasiconvex subgroup
of infinite index in $\Gamma$. Then the collection of translates of the (join of the) limit set of $H$ in $\partial \Gamma$ gives rise
to a symmetric pattern.\\

{\it 
For the rest of this subsection $Q$ will denote a group of {\bf boundary values} of pattern-preserving quasi-isometries of (the Cayley graph of) $G$.}\\

Recall that the topology on $Q$ is inherited from
 the uniform topology on $Homeo (\partial \Gamma )$. Also recall that $q \in Q$ is a 
quasi-symmetric map on the 
boundary $\partial \Gamma$. 

\begin{lemma} $q_n \rightarrow q$ in the uniform topology on $Homeo (\partial \Gamma )$
if and only if  $q_n^3 \rightarrow q^3$  in the compact-open topology on $\partial^3 \Gamma$. \end{lemma}

\begin{proof} Suppose   $q_n \rightarrow q$  in the uniform topology on $Homeo (\partial \Gamma )$. Then 
$q_n^3 \rightarrow q^3$  in the uniform topology on $Homeo (\partial \Gamma \times \partial \Gamma \times \partial \Gamma)$.
Since $\partial^3 \Gamma$ is an open invariant subset of $\partial \Gamma \times \partial \Gamma \times \partial \Gamma$, it follows that
 $q_n^3 \rightarrow q^3$  on compact subsets of $\partial^3 \Gamma$. Hence 
 $q_n^3 \rightarrow q^3$  in the compact-open topology on $\partial^3 \Gamma$.

Next suppose  $q_n^3 \rightarrow q^3$  in the compact-open topology on $\partial^3 \Gamma$. 
If $q_n$ does not converge to $q$ uniformly on $\partial \Gamma$, then (after passing to a subsequence if necessary)
there exists $\epsilon > 0$ such that for all $n$ there exists
$x_n \in \partial \Gamma$ such that $d(q_n(x_n) , q(x_n)) \geq 2\epsilon$. Passing to a further subsequence if necessary we can assume that 
$x_n \rightarrow x \in \partial \Gamma$. Hence $q(x_n) \rightarrow q(x) \in \partial \Gamma$. Hence by passing to a 
 further subsequence if necessary we can assume that $d(q_n(x_n) , q(x)) \geq \epsilon$ for all $n$. Choosing $y, z$ fixed distinct points unequal to $x$, it follows that 
$q_n^3(x_n,y,z)$ does not converge to $q^3(x,y,z)$, a contradiction.
\end{proof}

\begin{prop} $Q$ has no non-trivial topologically infinitely divisible elements. More generally, $Q$ does not
contain any non-trivial infinitely divisible subgroups.
\label{infdiv}
\end{prop}

\noindent {\bf Proof:} Suppose $q \in Q$ is infinitely divisible. Then for every neighborhood $U$ of $1$ there exists 
 $m \in \natls$ such that $q \in U^{m}$. Now for any finite collection $L_1, \cdots , L_n \in \LL$, there exists a neighborhood 
$U$ of $1 \in Q$ such that if  $g \in U$ then $g (L_i) = L_i$ for $i = 1, \cdots , n$. 

Hence for any finite collection $L_1, \cdots , L_n \in \LL$,
there exists a neighborhood $U$ of $1$
and $m \in \natls$ such that\\
1) $q \in U^{m} $ and \\
2) $g (L_j) = L_j$ for $g \in U, j = 1, \cdots , n$. \\

From (2) it follows that $g (L_j) = L_j$ for $g \in U^m$, $j = 1, \cdots , n$ and all $m \in \mathbf{N}$. \\
Therefore $q(L_j) =  L_j$ for all
$L_j \in \LL$. That is $q$ stabilizes every $L \in \LL$. If $x \in \partial \Gamma$, there exist $L_m \in \LL$ such that $L_m \rightarrow \{ x \}$
(the singleton set containing $x$) in the Hausdorff topology on $\bdy \Gamma$. Therefore $q( \{ x \} ) = \{ x \}$ for all $x \in \partial \Gamma$,
i.e. $q$ is the trivial element of $Q$. The same argument shows that $Q$ has no non-trivial topologically infinitely divisible subgroups. $\Box$

\begin{rmk} We state the second conclusion of Proposition \ref{infdiv} slightly differently. Let $U_i$ be a decreasing
sequence of symmetric neighborhoods of the identity in $Q$ such that $\bigcap_i U_i = \{ 1 \}$. Let
$\langle U_i \rangle = \bigcup_n U_i^n$. Then $\bigcap_i \langle U_i \rangle = \{ 1 \}$.
\end{rmk}

Since a connected topological group is generated by any neighborhood of the identity, we obtain

\begin{prop} $Q$ is totally disconnected. 
\label{td}
\end{prop}

Note that in Proposition \ref{infdiv}
and Proposition \ref{td} we {\bf do not} need to assume that $Q$ is a group of {\bf uniform} quasi-isometries.

\begin{rmk} \label{permutn} {\rm Let $K$ be a compact 
group of pattern-preserving quasi-isometries.
Then a reasonably explicit structure of $K$ may be given  as a permutation group. Since $K$ is compact
it acts on the discrete set $\LL$ with compact and hence finite orbits. Let $\LL_1, \LL_2, \cdots$ be a decomposition
of $\LL$ into disjoint orbits under $K$. Then $K \subset \Pi_i S(\LL_i)$, where $S(\LL_i )$ denotes the symmetric
group on the finite set $\LL_i$ and $\Pi$ denotes direct product. Thus, we have a natural representation of
$K$ as a permutation group on an infinite set, where every orbit is finite. The last part of the argument in
Proposition \ref{infdiv} shows that this representation is faithful, since any element stabilizing
every element of $\LL$ must be the identity.}
\end{rmk}

\subsection{PD Groups}
Boundaries $\partial G$
of PD(n) hyperbolic
groups $G$ are locally connected
homological manifolds (over the integers) with the homology
of a sphere of dimension $(n-1)$ by Theorem \ref{bestvina}. The interested reader may refer to Davis' survey \cite{davis-pd} for background
 on PD groups.

If $QI(G)$ denotes the group of boundary values
of quasi-isometries of $G$ acting on $\partial G$ (equipped with the uniform topology), then Theorem \ref{ns} implies
the following.

\begin{prop}
$QI(G)$ cannot have arbitrarily small torsion elements.
\label{pt}
\end{prop}

We combine Corollary \ref{hscor} with Proposition \ref{pt} to get the following.
\begin{prop} 
Let $G$ be a Poincare duality hyperbolic group and and $H$ an  infinite quasiconvex subgroup
of infinite index in $G$. Let $K$ be a compact group of (boundary values of) pattern-preserving
quasi-isometries of $G$. Suppose $d$ is a visual metric on $\partial G$ with 
$ dim_{haus} < dim_t +2$, where $dim_{haus}$ is the
Hausdorff dimension and $dim_t$ is the
topological dimension of $(\partial G, d)$.  Equivalently suppose that $ACD( \partial G) < dim_t + 2$ where $ACD(\partial G)$ denotes the Ahlfors regular conformal dimension of $\partial G$.  Then  $K$
 must be finite.
\label{cpt=finite}
\end{prop}

\noindent {\bf Proof:}  By Proposition \ref{pt}, $K$ cannot have arbitrarily small torsion elements. Hence there exists
$\epsilon > 0$ such that for any non-trivial $k \in K$, there exists an $x \in \partial G$ such that $\langle k \rangle x$ has diameter greater
than $\epsilon$ where  $\langle k \rangle$ denotes the cyclic group generated by $k$. Again, since $K$ is infinite and compact
there exists a sequence of distinct elements $k_i \rightarrow 1$ in $K$. By compactness of $\partial G$, there exists $x \in \partial G$,
such that $\langle k_i \rangle x$ has diameter greater
than $\epsilon$ and hence the order $o(k_i) \rightarrow \infty$ as $i \rightarrow \infty$. 

Therefore the subgroups  $\langle k_i \rangle \subset K$ converge up to a subsequence
(in the Chabauty topology on closed subgroups of $K$) to an infinite compact non-trivial abelian group $L$ without 
small torsion elements. If $L$ is pure torsion it must have elements of arbitrarily large order
(since $L$ being infinite and compact must contain elements arbitrarily close to the identity). This contradicts the structure
of compact abelian torsion groups (see Theorem 25.9 of \cite{hewitt-ross1}).

Hence, if $K$ is infinite, it must have an element $g$ of infinite order. Let $C(g)$ be the (closed) monothetic
subgroup generated by $g$. Since $K$ is totally disconnected by Proposition \ref{td} so is $C(g)$
and hence $C(g)$ cannot have arbitrarily small torsion elements. 
By Theorem \ref{div} $C(g)$ must contain a copy of the $p$-adic integers.
But $K$ cannot contain a copy of the $p$-adic integers by Corollary \ref{hscor}, a contradiction.
Hence $K$  is finite. $\Box$

We come now to the main
 Theorem of this section. Since $G$ acts on its Cayley graph $\Gamma$ by isometries, we are interested in uniform pattern-preserving groups 
of quasi-isometries containing $G$.

\begin{theorem} Let $G$ be a hyperbolic Poincare duality group and $H$ an  infinite quasiconvex subgroup
of infinite index in $G$. Suppose $d$ is a visual metric on $\partial G$ with 
$ dim_{haus} < dim_t +2$, where $dim_{haus}$ is the
Hausdorff dimension and $dim_t$ is the
topological dimension of $(\partial G,d)$.
Equivalently suppose that $ACD(\partial G) 
< dim_t + 2$ where $ACD(\partial G)$ denotes the Ahlfors regular conformal dimension of $\partial G$.
Let $Q$ be a group 
of pattern-preserving uniform quasi-isometries containing $G$. Then $G$ is of finite index in $Q$. In particular,
$Q \subset Comm (G)$, where $Comm (G)$ denotes the abstract commensurator of $G$.
\label{weakrig}
\end{theorem}

\noindent {\bf Proof:} 
Let $L$ be the limit set of $H$ and $\LL$ be the collection of translates of $L$ under $G$. By Corollary
\ref{upp-cpt}, we can choose a finite collection  $L_1 \cdots L_n$ 
  of elements of $\LL$ such that $Q_0 = \cap_{i = 1 \cdots n} Stab(L_i)$ is compact, where $Stab(L_i)$ denotes the
stabilizer of $L_i$ in $Q$. Then $Q_0$ is finite by Proposition \ref{cpt=finite}. As in the proof of Corollary
\ref{upp-lc}, we can choose a neighborhood $U$ of the identity in $Q$ such that $U \subset Q_0$. Hence $U$ is finite
and $Q$ is discrete.

Let $Gq_1, \cdots Gq_n, \cdots$ be distinct cosets. Since $G$ acts transitively on (the vertex set of) $\Gamma$, we can choose
representatives $g_1q_1, \cdots g_nq_n, \cdots$ such that $g_iq_i (1) = 1$ for all $i$. Since (the vertex set of) $\Gamma$ is locally finite,
the sequence $g_1q_1, \cdots g_nq_n, \cdots$ must have a convergent subsequence in $Q$. Since $Q$ is discrete, it follows that such
a sequence must be finite. Hence $G$ is of finite index in $Q$. 

Let $q \in Q$. Since $G$ is of finite index in $Q$, it follows that $qGq^{-1}$ is of finite index in $qQq^{-1} = Q$.  Therefore
$G \cap qGq^{-1}$ is of finite index in $Q$. Hence  $G \cap qGq^{-1}$ is of finite index in both $G$
and $qGq^{-1}$, i.e. $q \in Comm (G)$, where $Comm (G)$ denotes the abstract commensurator of $G$ (where we identify $q$
with the element of $Comm(G)$ that takes  $G \cap qGq^{-1}$ to $q(G \cap qGq^{-1})q^{-1}$). Since distinct elements $q_1, q_2$ induce
distinct homeomorphisms of $\partial G$ by definition, and since two elements defining the same element of
$Comm (G)$ induce the same homeomorphism on $\partial G$, the elements $q_1, q_2 \in Comm (G)$ are distinct. This proves the result. $\Box$

\smallskip

In fact the proof of Theorem \ref{weakrig} gives:

\begin{cor}
Let $G$ be a hyperbolic Poincare duality group and $H$ an  infinite quasiconvex subgroup
of infinite index in $G$. Suppose $d$ is a visual metric on $\partial G$ with 
$dim_{haus} < dim_t +2$, where $dim_{haus}$ is the
Hausdorff dimension and $dim_t$ is the
topological dimension of $(\partial G,d)$.
Equivalently suppose that $ACD(\partial G) 
< dim_t + 2$ where $ACD(\partial G)$ denotes the Ahlfors regular conformal dimension of $\partial G$.
Let $Q$ be a  locally compact group 
of pattern-preserving quasi-isometries containing $G \subset Homeo(\partial G)$, where $Q$ is 
equipped with the uniform topology. Then $G$ is of finite index in $Q$. In particular,
$Q \subset Comm (G)$, where $Comm (G)$ denotes the abstract commensurator of $G$.
\label{weakriglc}
\end{cor}

\noindent {\bf Proof:} By Proposition \ref{td}, $Q$ is totally disconnected. Hence by Theorem \ref{hr1}
$Q$ contains arbitrarily small compact open subgroups. Such compact open subgroups are
 finite by Proposition \ref{cpt=finite} (as in the proof of Theorem \ref{weakrig}). Hence $Q$ is discrete.
The rest of the proof is as in Theorem \ref{weakrig}. $\Box$

\smallskip

In the present context, Scholium \ref{gromov} translates to the following precise statement as a consequence
of Corollary \ref{weakriglc}.

\begin{cor} Let $\phi$ be a pattern-preserving quasi-isometry between pairs $(G_1, H_1)$ and $(G_2, H_2)$
of hyperbolic PD groups and infinite quasiconvex subgroups of infinite index. Suppose $d$ is a visual metric on 
$\partial G_1$ with 
$dim_{haus} < dim_t +2$, where $dim_{haus}$ is the
Hausdorff dimension and $dim_t$ is the
topological dimension of $(\partial G_1,d)$. 
Equivalently suppose that $ACD(\partial G_1) 
< dim_t + 2$ where $ACD(\partial G_1)$ denotes the Ahlfors regular conformal dimension of $ (\partial G_1)$.
Further, suppose that $G_1$
and $\partial \phi^{-1} \circ G_2 \circ \partial \phi$ embed in  some locally compact subgroup $Q$ of
$Homeo (\partial G_1)$ with the uniform topology. Then $G_1$ and $G_2$ are commensurable.
\label{precisescholium}
\end{cor}

\noindent {\bf Proof:} By Corollary \ref{weakriglc}, $G_1$ is of finite index in $Q$. Hence $Q$ is a rational PD(n) group. Since  
 $\partial \phi^{-1} \circ G_2 \circ \partial \phi$ embeds in $Q$ and is also a PD(n) group, it follows that 
$\partial \phi^{-1} \circ G_2 \circ \partial \phi \subset Q$ is of finite index \cite{brown}. Hence 
$G_1$ and $G_2$ are commensurable.
$\Box$

\section{Filling Codimension One Subgroups and Pattern Rigidity}

\subsection{Codimension One Subgroups and Pseudometrics}\label{codimone}
Let $G$ be a one-ended Gromov-hyperbolic group with Cayley graph $\Gamma$. 
Let $H$ be a quasiconvex subgroup. We say that $H$ is {\bf codimension one} if the limit set $L_H$
of $H$ disconnects $\bdy G$. This is equivalent to saying that the join $J(L_H) = J$ disconnects $\Gamma$
coarsely, i.e. if $D$ be the quasiconvexity constant of $J$, then
$\Gamma \setminus N_D(J)$ has more than one unbounded component, where $N_D(J)$ denotes the $D$-neighborhood
of $J$.  (See Ch. 2 \cite{ss-book}, particularly Remark 2.4 for a proof of this equivalence and related results.)

We say further that $H$ (or more generally a finite collection $H_1, \cdots H_k$)
is {\bf filling} if for any two $x , y \in \bdy G$, there exists a translate $gL_H$ of
$L_H$ (or more generally $ gL_{H_i}$) by an element $g$ of $G$ such that $x, y$ lie in distinct components of $\bdy G \setminus gL_H$
(or more generally $\bdy G \setminus gL_{H_i}$). We shall deal with a single filling subgroup for convenience
of exposition. The results in this section go through for a finite collection of filling codimension one subgroups.

The existence
of a finite collection of filling codimension one quasiconvex subgroups is important in light of the following Theorem
due to Sageev \cite{sageev} and Bergeron-Wise \cite{berg-wise}.

\begin{theorem} \cite{sageev} \cite{berg-wise} A hyperbolic group
acts properly, cocompactly on a CAT(0) cube complex if and only if it admits a  a finite collection $H_1, \cdots H_k$ of 
filling codimension one quasiconvex subgroups. \label{sbw} \end{theorem}

Let $D_1$ be such that any path joining points in distinct unbounded components of 
$\Gamma \setminus N_D(J)$ passes within $D_1$ of $J$. 
We say that $x , y \in \Gamma$ are separated by some translate $gJ$ of $J$ if $x, y$ lie in distinct unbounded
components of $\Gamma \setminus gN_{D+D_1}(J)$. Equivalently, we shall say that the geodesic $[x,y]$ is 
separated by some translate $gJ$ of $J$.

\begin{lemma} Let $H$ be a codimension one, filling, quasiconvex subgroup of a one-ended hyperbolic group $G$. 
Let $\Gamma$ be a Cayley graph of $G$. There exists 
$C\geq 0$ such that any geodesic $\sigma$ in $\Gamma$ of length greater than $C$ is separated by a translate of $J$.
\label{separate}
\end{lemma}

\noindent {\bf Proof:} Suppose not. Then there exists a sequence of geodesic segments $\sigma_i = [a_i, b_i]$ which are not separated by any
translate of $J$ such
that $d(a_i,b_i) \rightarrow \infty$.
By equivariance, we may assume that $\sigma_i$ is centered at the origin, i.e.  $d(a_i,1) \rightarrow \infty$ and
 $d(1,b_i) \rightarrow \infty$ and $1\in [a_i,b_i]$. Let $a_i \rightarrow a_\infty \in \bdy G$
and $b_i \rightarrow b_\infty \in \bdy G$. Then $a_\infty$ and $b_\infty$ cannot lie in distinct components of 
$\bdy G \setminus gL_H$ for any $g \in G$, for if they did then there exists $g \in G$ such that
$a_\infty$ and $b_\infty$  lie in distinct components of 
$\bdy G \setminus gL_H$ and hence
for all $i$ sufficiently large,
$a_i, b_i$ would lie in distinct unbounded components of $\Gamma \setminus gN_D(J)$. (Here we are implicitly using the correspondence
between the  unbounded components of $\Gamma \setminus gN_D(J)$ and the components of 
$\bdy G \setminus gL_H$ mentioned in the first paragraph of this subsection cf. Ch. 2 of \cite{ss-book}). 

But if 
$a_\infty$ and $b_\infty$ cannot be separated, then $H$ cannot be filling,
contradicting the hypothesis. $\Box$

\begin{lemma} Let $G, H, \Gamma, J$ be as above.
Let $[a,b] \subset \Gamma$ be a geodesic and $c \in [a,b]$ such that $d(a,c) \geq 2D$, 
$d(b,c) \geq 2D$, where $D$ is the quasiconvexity constant of $J$. Suppose $gJ$ separates $a, c$. Then 
$gJ$ separates $a, b$.
\label{monotonic} \end{lemma}

\noindent {\bf Proof:} Suppose not. Then $a, b$ lie in the same unbounded component of $\Gamma \setminus
gN_D(J)$, whereas $c$ lies in a different unbounded component of $\Gamma \setminus
gN_D(J)$. Hence there is a subsegment $[ecf]$ of $[a,b]$ such that $e, f \in N_{D_1}J$, but
$c \notin  N_{D+D_1}J$, contradicting the quasiconvexity constant for $J$. $\Box$

\begin{lemma} \cite{GMRS} Let $G$ be a hyperbolic group and $H$ a quasiconvex subgroup, with limit set
$L$. Let $J$ denote the join of the limit set. There exists $N \in \natls$ such that
there exist at most $N$ distinct translates of $J$ intersecting the $2$-neighborhood $B_2(g)$ nontrivially
for any $g \in G$.
\label{bdd}
\end{lemma}

\begin{lemma} Define a new pseudometric $\rho_1$ on $\Gamma$ by declaring  $\rho_1 (a,b)$ to be the number of copies of 
joins $J \in \JJ$ separating $a,b$. Then $(\Gamma, \rho_1 )$ is quasi-isometric to $(\Gamma , d)$
\label{qi}
\end{lemma}

\noindent {\bf Proof:}
By Lemma \ref{bdd}, 
it follows that  there exists $N \in \natls$ such that $d(a,b) \leq C_0$
implies $\rho (a,b) \leq NC_0$. From Lemma \ref{separate}, it follows
 that there exists $C_2 \geq 0$,  such that  $d(a,b) \geq C_2 $
implies $\rho (a,b)) \geq 1$. Now from Lemma \ref{monotonic}, it follows that for
$n \in \natls$, $d(a,b) \geq nC_2 $  implies $\rho (a,b)) \geq n$. Hence the Lemma. $\Box$

\smallskip

A purely topological version of Lemma \ref{qi} may be obtained as follows. Let $\partial^3G$ denote the collection
of distinct unordered triples of points on $\partial G$. Then it is well known \cite{gromov-hypgps} \cite{bowditch-jams}
that $G$ acts cocompactly on $\partial^3G$ with metrizable quotient. Let $\rho$ be the pseudo-metric of Observation \ref{milnor-s} which asserts
that $(\partial^3G, \rho )$ is quasi-isometric to $(\Gamma , d)$.

We say that a translate $gL \in \LL$ separates closed subsets $A, B \subset \partial G$ if $A,B$ lie in distinct
components of $\partial G \setminus gL$.
Define a pseudometric $\rho_2$ on $\partial^3G$ by defining $\rho_2 (\{a_1, a_2, a_3\},\{b_1, b_2, b_3\})$ to be the number of copies of 
limit sets $gL \in \LL$ separating $\{a_1, a_2, a_3\},\{b_1, b_2, b_3\}$. Then $(\partial^3G, \rho_2)$ is quasi-isometric
to $(\Gamma, \rho_1 )$, and hence to $(\Gamma, d )$ and $(\partial^3G, \rho)$. We state this as follows.

\begin{cor} $(\partial^3G, \rho_2)$,
 $(\Gamma, \rho_1 )$,  $(\Gamma, d )$ and $(\partial^3G, \rho)$ are  quasi-isometric to each other.
\label{qicor}
\end{cor}

\subsection{Pattern Rigidity}

We prove the following Proposition for which $G$ may be {\it any} one-ended hyperbolic group (not necessarily PD):

\begin{prop}
Let $G$ be a one-ended hyperbolic group and $H$ a codimension one, filling, quasiconvex subgroup. 
Then any pattern-preserving group $Q$
of quasi-isometries is uniform. 
\label{filling}
\end{prop}

\noindent {\bf Proof:}  Assume without loss of generality that $G \subset Q$. Since $Q$ consists of 
pattern-preserving quasi-isometries, each element
of $Q$ induces a pattern-preserving homeomorphism of $\partial G$.
Since any pattern-preserving homeomorphism of $\partial G$ preserves $(\partial^3G, \rho_2)$ on the nose, it follows
from Corollary \ref{qicor} that $Q$ is uniform. $\Box$

Combining Proposition \ref{filling} with Theorem \ref{weakrig} we  get the following Theorem. (Note that the only
PD group that is not one-ended is $\Z$, in which case codimension one quasiconvex subgroups in our sense do not exist.)

 \begin{theorem}
Let $G$ be a PD hyperbolic group and $H$ a codimension one, filling, quasiconvex subgroup. 
Let  $Q$ be any pattern-preserving group
of quasi-isometries containing $G$. Suppose $d$ is a visual metric on 
$\partial G$ with 
$dim_{haus} < dim_t +2$, where $dim_{haus}$ is the
Hausdorff dimension and $dim_t$ is the
topological dimension of $(\partial G,d)$. 
Equivalently suppose that $ACD(\partial G) 
< dim_t + 2$ where $ACD(\partial G)$ denotes the Ahlfors regular conformal dimension of $\partial G$.
Then  the index of $G$ in $Q$ is finite.
\label{strongrig}
\end{theorem}

In fact more is true. Combining Proposition \ref{filling} with Corollary \ref{qicor}, we get

\begin{prop}
Let $G$ be a one-ended hyperbolic group and $H$ a codimension one, filling, quasiconvex subgroup with limit set
$L$. Let $\LL$ be the collection of translates of $L$ under $G$. 
Then any pattern-preserving group $Q$
of homeomorphisms of $\partial G$  preserving $\LL$ can be realized as the boundary values
of uniform quasi-isometries. 
\label{filling2}
\end{prop}

Note that in Proposition \ref{filling2} we do not need $G$ to be a PD group.
Combining Proposition \ref{filling2} with Theorem \ref{weakrig} we finally get

 \begin{theorem} {\bf Topological Pattern Rigidity}
Let $G$ be a PD hyperbolic group and $H$ a codimension one, filling, quasiconvex subgroup
 with limit set
$L$. Let $\LL$ be the collection of translates of $L$ under $G$.  Suppose $d$ is a visual metric on 
$\partial G$ with 
$dim_{haus} < dim_t +2$, where $dim_{haus}$ is the
Hausdorff dimension and $dim_t$ is the
topological dimension of $(\partial G,d)$.  Equivalently suppose that $ACD(\partial G) 
< dim_t + 2$ where $ACD(\partial G)$ denotes the Ahlfors regular conformal dimension of $\partial G$.
Let  $Q$ be any pattern-preserving group
of homeomorphisms of $\partial G$  preserving $\LL$ and containing $G$. Then  the index of $G$ in $Q$ is finite.
\label{toprig}
\end{theorem}

Theorem \ref{toprig} is a generalization of a Theorem of Casson-Bleiler \cite{CB} and Kapovich-Kleiner
\cite{kap-kl-lowd} to all dimensions. Casson-Bleiler \cite{CB} and Kapovich-Kleiner
\cite{kap-kl-lowd} proved Theorem \ref{toprig} for $G$ the fundamental group of a surface
and $H$ an infinite cyclic subgroup corresponding to a filling curve.

\section{Finite Volume Manifolds of Negative Curvature}

Let $M = M^n$ be a complete finite volume non-compact manifold of pinched negative curvature (i.e. $-1 \leq \chi \leq -K$ for some $K \geq 1$,
where $\chi$ denotes sectional curvature). Then $\widetilde{M}$ is homeomorphic to $\R^n$ by the Cartan-Hadamard Theorem and its
ideal boundary $\partial \widetilde{M}$ is homeomorphic to $S^{n-1}$. Further, by Theorem \ref{ahlfors2}
$X=\partial\widetilde{M}$ equipped with a visual metric has the structure of an Ahlfors regular metric measure space.
In fact, if $-1 \leq \chi \leq -(1+\epsilon )^2$, 
 the Hausdorff dimension of the visual boundary is bounded above by $(n-1)(1+\epsilon )$ (see Remark \ref{pinch} below).

\subsection{Symmetric Patterns of Horoballs} $M$ has a finite number of cusps. Lifting these to $\widetilde{M}$ we obtain an
equivariant collection $c\HH$ of horoballs. For convenience of exposition we assume that $M$ has one cusp. 
We shall denote   individual
elements of $c\HH$ by $cH$ or $cH_i$. The boundary of the horoball $cH$ is called a horosphere 
and is denoted as $H$. The collection of horospheres will be denoted by $\HH$.
Let $G = \pi_1(M)$ and let $K $ denote the fundamental group of the cusp.

The collection $c\HH$ will be called a symmetric pattern
(of horoballs).
It is a fact that elements of $c\HH$ are uniformly quasiconvex \cite{farb-relhyp} and that for any two
distinct $cH_1, cH_2 \in c\HH$, there is a coarsely well-defined `centroid', i.e. the shortest
geodesic joining $cH_1, cH_2 \in c\HH$ is  coarsely well-defined (any two such lie in a uniformly
bounded neighborhood of each other \cite{farb-relhyp} )
and hence its mid-point (the centroid of $cH_1, cH_2$) is coarsely well-defined.

In \cite{farb-relhyp}, Farb proves that $G$ is strongly hyperbolic relative to $K$. This is equivalent to the statement  that
$\widetilde{M}$ is strongly hyperbolic relative to $c\HH$. Equivalently, the `neutered space' 
$\widetilde{M} \setminus ~ \bigcup_{cH \in c\HH} ~ Int ~(cH)$ is hyperbolic relative to the collection $\HH$.
We refer to  \cite{farb-relhyp} for background on relative hyperbolicity.

We now recast the relevant
definitions and propositions of Sections 1 and 3 in the present context.
Let $\Gamma$, $\Gamma_K$, $\Gamma_\KK$ denote respectively the Cayley graph of $G$,
some translate of the Cayley (sub)graph of $K$ and the collection of translates of $\Gamma_K$ (assuming as usual that the finite generating set
of $K$ used in constructing $\Gamma_K$ is contained in the the finite generating set
of $G$ used in constructing $\Gamma$).

\begin{defn} The group $PP(G,K)$ of pattern-preserving maps for a (strongly) relatively
hyperbolic pair $(G,K )$ as above is defined as the group
of homeomorphisms of $X=\partial\widetilde{M}$  preserving (as a set) the collection of base-points
of $c\HH$. 
The group $PPQI(G,K)$ of pattern-preserving quasi-isometries
 for such a (strongly) relatively
hyperbolic pair $(G,K )$  is defined as the subgroup of $PP(G,K)$ consisting of homeomorphisms $h$
of $\partial G$ such that $h = \partial \phi$ for some quasi-isometry $\phi : 
\widetilde{M} \rightarrow \widetilde{M}$ that permutes the collection
of horoballs $c\HH$.
\end{defn}

The following Theorem is a special case of a Theorem  proven in \cite{mahan-relrig} using the notion of {\it mutual coboundedness}.

\begin{theorem}\cite{mahan-relrig} Let $M = M^n$ be a complete finite volume  manifold of pinched negative curvature and let
$c\HH$ denote the associated symmetric pattern of horoballs.

 There exist two elements
$cH_1,  cH_2$ of  $c\HH$  such that the following holds. \\
For any $K, \epsilon$, there exists a $C$ such that if $\phi : \widetilde{M} \rightarrow \widetilde{M}$ is a pattern-preserving
$(K, \epsilon )$-quasi-isometry  with $\partial \phi (\partial cH_i) = \partial cH_i$ for $i = 1, 2$, then $d(\phi (1), 1) \leq C$.
\label{originrh}
\end{theorem}

Let $\overline{M} = \widetilde{M} \cup \partial \widetilde{M}$ denote the Gromov compactification of $\widetilde{M}$
and $\overline {c\HH}$ denote the collection of compactified horoballs, i.e. horoballs with basepoints adjoined.
Let $d_c$ denote a metric giving the topology on $\overline{M}$.
In this context Proposition \ref{discrete} translates to the following (see \cite{mahan-relrig} for instance).

\begin{prop} 
The collection $\overline {c\HH}$ is discrete in the Hausdorff topology on the space of closed subsets of 
$\overline{M}$, i.e. for all $cH \in c\HH$,
there exists $\epsilon > 0$ such that $N_\epsilon (cH) \cap c\HH = cH$, where  $N_\epsilon (cH) $ denotes an $\epsilon$ neighborhood of $cH$
in the Hausdorff metric arising from $d_c$. \label{discreterh} \end{prop}

Let $Q \subset PPQI(G,K)$ be a group of quasi-isometries preserving a symmetric pattern of horoballs.
Using Theorem \ref{originrh} and Proposition \ref{discreterh}, we have as in Section 3 (cf Proposition \ref{infdiv},
Corollary \ref{td}, Proposition \ref{pt}, Proposition \ref{cpt=finite}):

\begin{prop} $Q$ has no non-trivial topologically infinitely divisible elements. More generally, $Q$ does not
contain any non-trivial infinitely divisible subgroups. Hence $Q$ is totally disconnected. 
Suppose  further that $\partial \widetilde{M}$ has a visual metric $d$  with 
$dim_{haus} < dim_t +2$, where $dim_{haus}$ is the
Hausdorff dimension and $dim_t$ is the
topological dimension of $(\partial \widetilde{M},d)$. Equivalently suppose that $ACD(\partial (\widetilde{M})) 
< dim_t + 2$ where $ACD(\partial (\widetilde{M}))$ denotes the Ahlfors regular conformal dimension of $\partial (\widetilde{M})$. If
$Q$ is compact, then $Q$
is finite.
\label{omnirh}
\end{prop}

\noindent {\bf Sketch of Proof:} 
Proposition \ref{infdiv} and
Corollary \ref{td}  apply directly to show that $Q$ does not have
 non-trivial topologically infinitely divisible elements.

Theorem \ref{ns} shows that $Q$ does not have arbitrarily small torsion elements.

The rest of the argument is as in
  Proposition \ref{cpt=finite}.
The only point that needs to be mentioned is that $\partial^3 \widetilde{M}$ minus a neighborhood of the cusps,
rather than  $\partial^3 \widetilde{M}$  itself
is  quasi-isometric to the Cayley graph of $\pi_1(M)$. Since the quasi-isometries in $Q$ preserve the horoballs they naturally
induce homeomorphisms of  $\partial^3 \widetilde{M}$
preserving horoballs. Observation \ref{milnor-s} therefore goes through with the above modification.
 $\Box$

\medskip

Now, let  $Q_u \subset PPQI(G,K)$ be a group of {\it uniform}
quasi-isometries preserving a symmetric pattern of horoballs. Then $Q_u$ is locally compact (by Lemma \ref{lc})
and contains a compact open subgroup $K_0$. 
Suppose  further that 
$d$ is a visual metric on 
$\partial (\widetilde{M} )$ with 
$dim_{haus} < dim_t +2$, where $dim_{haus}$ is the
Hausdorff dimension and $dim_t$ is the
topological dimension of $(\partial (\widetilde{M} ),d)$. 
Equivalently suppose that $ACD(\partial (\widetilde{M})) 
< dim_t + 2$ where $ACD(\partial (\widetilde{M}))$ denotes the Ahlfors regular conformal dimension of $\partial (\widetilde{M})$.
Then $K_0$ is finite by Proposition \ref{omnirh}. Hence $Q_u$ is discrete. Thus as in Theorem \ref{weakrig}
we get the following.

\begin{theorem} Let $M = M^n$ be a complete finite volume  manifold of pinched negative curvature and let
$c\HH$ denote the associated symmetric pattern of horoballs.
Suppose  further that 
$d$ is a visual metric on $\partial (\widetilde{M})$ 
 with 
$dim_{haus} < dim_t +2$, where $dim_{haus}$ is the
Hausdorff dimension and $dim_t$ is the
topological dimension of $(\partial  (\widetilde{M}),d)$. 
Equivalently suppose that $ACD(\partial (\widetilde{M})) < dim_t + 2$ where $ACD(\partial (\widetilde{M}))$ 
denotes the Ahlfors regular conformal dimension of $\partial (\widetilde{M})$. Let $Q$ be a group 
of  uniform quasi-isometries containing $G$ preserving a symmetric pattern of horoballs. Then $G$ is of finite index in $Q$. In particular,
$Q \subset Comm (G)$, where $Comm (G)$ denotes the abstract commensurator of $G$.
\label{weakrigrh}
\end{theorem}

\subsection{Weak QI rigidity for Relatively Hyperbolic Groups}

We shall be using the following Theorem of  Behrstock-Drutu-Mosher (which follows from the
proof of Theorem 4.8 of \cite{BDM}; see also \cite{schwartz_pihes}) and Farb's result \cite{farb-relhyp} that the fundamental group
of complete finite volume non-compact manifold of pinched negative curvature is strongly hyperbolic relative to the cusp groups.

\begin{theorem} \label{bdm} {\bf (Behrstock-Drutu-Mosher \cite{BDM}, Schwartz \cite{schwartz_pihes})} 
Let $M = M^n$ be a complete finite volume  manifold of pinched negative curvature with $n > 2$ (and one cusp
for ease of exposition). Let $G= \pi_1(M)$.
Let $K$ denote the fundamental group of the cusp. Choose a finite generating set for $G$ containing a 
finite generating set for $K$.
 Let $\Gamma$, $\Gamma_K$, $\Gamma_\KK$ denote respectively the Cayley graph of $G$,
some translate of the Cayley (sub)graph of $K$ and the collection of translates of $\Gamma_K$. Then for every
$L \geq 1$ and $C \geq  0$ there exists $R = R(L, C, G, K )$ such that the following holds. \\
For any $(L, C)$- (self) quasi-isometry $q$ of $G$,
the image $q(\Gamma_K)$ is at a bounded Hausdorff distance $R$ of some $g\Gamma_K  \in \Gamma_\KK$.
\end{theorem} 

Let $q$ be a (self) quasi-isometry
of $\Gamma$. Elements of
the collection $\Gamma_\KK$ are mapped bijectively to bounded neighborhoods
of elements of
 collection $\Gamma_\KK$. Identifying $\Gamma$ with the neutered space 
$(\widetilde{M} \setminus \bigcup_{cH \in c\HH} ~ Int (cH))$ we can extend
 $q$  (cf. \cite{schwartz_pihes}) to a (self) quasi-isometry $q^h$ of $\widetilde{M}$
where the elements of $c\HH$ are bijectively mapped to uniformly bounded neighborhoods of elements of $c\HH$.
Each element of $c\HH$ has a unique limit point in $\partial \widetilde{M}$ which we shall call its base-point.
Let $\partial q$ denote the induced map of $\partial (\widetilde{M} )$ and $\partial \HH$ denote the collection
of base-points of $c\HH$ in $\partial \widetilde{M}$.
Thus 
 a simple consequence
of Theorem \ref{bdm} is the following.

\begin{cor} \label{bdmcor}  Let $M = M^n$ be a complete finite volume  manifold of pinched negative curvature with $n > 2$ (and one cusp
for ease of exposition). Let $G= \pi_1(M)$.
Let $K$ denote the fundamental group of the cusp. Let $c\HH$ be the associated symmetric pattern of
horoballs in $\widetilde{M}$. Choose a finite generating set for $G$ containing a 
finite generating set for $K$.
 Let $\Gamma$, $\Gamma_K$, $\Gamma_\KK$ denote respectively the Cayley graph of $G$,
some translate of the Cayley (sub)graph of $K$ and the collection of translates of $\Gamma_K$. 
Identify $\Gamma$ (coarsely) with the neutered space 
$(\widetilde{M} \setminus \bigcup_{cH \in c\HH} ~ Int (cH))$.

Then for every
$L \geq 1$ and $C \geq  0$ there exist $L_1 \geq 1$, $C_1 \geq  0$ and $R = R(L,C,G, K )$
such that the following holds. \\
For any $(L, C)$- (self) quasi-isometry $q$ of $\Gamma$ (or equivalently, $(\widetilde{M} \setminus \bigcup_{cH \in c\HH} ~ Int (cH))$), 
there is an $(L_1, C_1)$ (self) quasi-isometry 
$q^h$ of $\widetilde{M}$ such that 
the image $q^h(\Gamma_K)$ is at a bounded Hausdorff distance $R$ of some $g\Gamma_K  \in \Gamma_\HH$.
Hence $q$ induces a homeomorphism $\partial q$ of $\partial (\widetilde{M} )$ preserving the base-points
of $c\HH$.
\end{cor} 

Combining Corollary \ref{bdmcor} with Theorem \ref{weakrigrh} we get the following.

\begin{theorem} Let $M = M^n$ be a complete finite volume  manifold of pinched negative curvature with $n > 2$. 
Let $G= \pi_1(M)$. Suppose that there exists
 a visual metric $d$ on 
$\partial (\widetilde{M})$ with 
$dim_{haus} < dim_t +2$, where $dim_{haus}$ is the
Hausdorff dimension and $dim_t$ is the
topological dimension of $(\partial (\widetilde{M}),d)$. 
Equivalently suppose that $ACD(\partial (\widetilde{M})) < dim_t + 2$ where $ACD(\partial (\widetilde{M}))$
 denotes the Ahlfors regular conformal dimension of $\partial (\widetilde{M})$. \\
Let $\Gamma$ be a Cayley graph of $G$ with respect to a finite generating set. Let $Q$ be a group 
of  uniform quasi-isometries of $\Gamma$ containing $G$. Then $G$ is of finite index in $Q$. In particular,
$Q \subset Comm (G)$, where $Comm (G)$ denotes the abstract commensurator of $G$.
\label{weakqirigrh}
\end{theorem}

\section{Examples and Consequences}

\subsection{Examples} In this subsection we list a collection of examples of hyperbolic $PD(n)$ groups.

\noindent {\bf Rank one Symmetric Spaces}\\
Uniform lattices in rank one symmetric spaces are examples of hyperbolic $PD(n)$ groups. For 
lattices in real hyperbolic space $\Hyp^n$ the usual visual metric has Hausdorff dimension
$dim_{haus}$ equal to topological dimension $dim_t$. For 
lattices in complex hyperbolic space $\C\Hyp^n$ the usual visual metric has 
$dim_{haus} = dim_t + 1$. 
Lattices in quaternionic hyperbolic space and the Cayley hyperbolic  plane cannot have codimension one 
subgroups; so Theorem \ref{weakrig} cannot apply. But these spaces are qi-rigid by a deep
Theorem of Pansu \cite{pansu-rig}.

\smallskip

\noindent {\bf Gromov-Thurston Examples}\\
In \cite{gt}, Gromov and Thurston construct examples of closed negatively curved $n$-manifolds ($n > 3$) of 
{\it arbitrarily} pinched negative
curvature, which do not admit metrics of constant negative curvature. 

\begin{rmk} \label{pinch} {\rm For these examples
the inequality
$dim_{haus} < dim_t + 2$ is satisfied. This is because of the following. If $-1 \leq \chi \leq -(1+\epsilon )^2$, the rate of divergence
of geodesics is bounded by $O(e^{(1+\epsilon )R})$. Hence volumes of $R$-balls is bounded by $O(e^{(n-1)(1+\epsilon )R})$.
Therefore the Hausdorff dimension of the visual boundary is bounded above by $(n-1)(1+\epsilon )$ (See Theorem 1 of
\cite{sullivan-acta}).} \end{rmk}

\smallskip

\noindent {\bf Mostow-Siu Examples}\\
In \cite{mostow-siu},  Mostow and Siu constructed an infinite family of complex surfaces that admit 
negatively curved Kahler metrics but do not admit a complex hyperbolic structure. Fundamental
groups of these provide further examples of hyperbolic $PD(n)$ groups.

\smallskip

\noindent {\bf Davis-Januskiewicz Examples}\\
    A remarkable set of examples is constructed by   
Davis and Januskiewicz \cite{dj}, who show that there exist hyperbolic
$PD(n)$ groups $G$ for $n \geq 4$, such that the boundary
$\partial G$ is not homeomorphic to the sphere $S^{ n-1}$.
 $\partial G$ need not be simply connected or locally
simply connected and hence is not even an Absolute Neighborhood Retract.

\subsection{Quasi-isometric Rigidity}

Let $\AAA$ be  a graph of groups  with Bass-Serre tree of spaces
$X \rightarrow T$. Let $A = \pi_1 \AAA$.    Let $\VV \EE (T )$
 be the set of vertices and edges of $T$. The metric on $T$ induces a metric
on $\VV \EE (T )$, via a natural injection $\VV \EE (T ) \rightarrow T$
 which takes each vertex to itself and each
edge to its midpoint. Let $d_H$ denote Hausdorff distance. 

We refer the reader to \cite{msw2} specifically for the following notions:\\
1) Depth zero raft. \\
2) Crossing graph condition. \\
3) Coarse finite type and coarse dimension.\\
4) Finite depth.

Combining Theorems
1.5, 1.6 of Mosher-Sageev-Whyte \cite{msw2}  with the Pattern Rigidity theorem
\ref{strongrig} we have the following QI-rigidity Theorem along the lines
of Theorem 7.1 of \cite{msw2}.\\

\begin{theorem} \label{qirig}
Let $\AAA$ be a finite, irreducible graph of groups such that for
the associated
Bass-Serre tree $T$ of spaces  \\
a) the vertex groups are $PD(n)$ hyperbolic groups for some fixed $n$ \\
b) edge groups are filling codimension one
in the adjacent vertex groups  \\
and such that $\AAA$ is of finite depth.
Further suppose that each vertex group $G$ admits  a visual metric $d$ on 
$\partial G$ with 
$dim_{haus} < dim_t +2$, where $dim_{haus}$ is the
Hausdorff dimension and $dim_t$ is the
topological dimension of $(\partial G,d)$. Equivalently suppose that $ACD(\partial G) < dim_t + 2$ where $ACD(\partial G)$ denotes the Ahlfors regular conformal dimension of $\partial G$.\\
 If $H$ is a finitely generated group quasi-isometric to $A = \pi_1 ( \AAA )$ then $H$ splits as
a graph $\AAA^{\prime}$ of groups whose depth zero vertex groups are commensurable to the 
depth zero vertex groups of $\AAA$ and
whose edge groups and positive depth vertex groups are  quasi-isometric to
groups of type   (b).
\end{theorem}

\noindent {\bf Proof:} By the restrictions on the vertex 
and edge groups, it automatically follows that
 all vertex and edge groups are PD groups of coarse finite type.
Since the edge groups are filling, the crossing graph condition  of Theorems 1.5, 1.6
of  \cite{msw2} is satisfied.
 $\AAA$ is
 automatically finite depth, because an infinite index subgroup of
a $PD(n)$ groups has coarse dimension at most $n - 1$. 

Then by  Theorems 1.5 and 1.6 of \cite{msw2}, 
$H$ splits as a graph of groups $\AAA^{\prime}$
 with depth zero vertex spaces
quasi-isometric to the vertex groups of $\AAA$ and edge groups quasi-isometric to 
the edge groups of $\AAA$. 
Further, the quasi-isometry respects the vertex and edge spaces of this 
splitting, and thus the quasi-actions of the vertex groups on the vertex spaces of $\AAA$ preserve
the patterns of edge spaces.

By Theorem \ref{strongrig} the depth zero vertex groups in  $\AAA^{\prime}$
 are commensurable to the corresponding groups in  $\AAA$. $\Box$

\subsection{The Permutation Topology}

In this paper we have ruled out three kinds of elements from the group $PPQI(G,H)$ of pattern-preserving
quasi-isometries under appropriate hypotheses on $G$:\\
a) Elements that admit arbitrarily small roots (topologically divisible elements)\\
b) Arbitrarily small torsion elements (essentially Theorem \ref{ns})\\
c) Elements with arbitrarily large powers close to the identity (no copies of the $p$-adics)\\

In a sense (a) and (c) are phenomena that are opposite to each other. In hindsight, the previous works
on pattern rigidity \cite{schwarz-inv} \cite{biswas-mj} \cite{bi} exploited (a) in the context of an
ambient Lie group which automatically rules out (b) and (c).

We have used the fact that the group is pattern-preserving in a rather weak sense, only to conclude that
the group we are interested in is totally disconnected. In fact Theorem \ref{weakrig} generalizes readily
to show that a locally compact totally disconnected group of quasi-isometries containing $G$ must be a
finite extension of $G$ under appropriate hypotheses on $G$.
 The crucial hypothesis is local compactness on $PPQI(G,H)$ which can be removed under hypotheses on $H$
as in Theorem \ref{strongrig}. We would like to remove the hypothesis of local compactness in more general
situations.

Remark \ref{permutn} gives a reasonably explicit structure of $K$ for $K$  a compact 
group of pattern-preserving quasi-isometries.
 $K$ 
 acts on the discrete set $\LL$ of patterns with  finite orbits $\LL_1, \LL_2, \cdots$ 
and hence $K \subset \Pi_i S(\LL_i)$, where $S(\LL_i )$ denotes the symmetric
group on the finite set $\LL_i$ and $\Pi$ denotes direct product. 

To establish Theorem \ref{weakrig} without the hypothesis of local compactness, two crucial problems remain:\\
{\bf Problem 1:} A topological converse to Remark \ref{permutn} which would say that a group
$K \subset \Pi_i S(\LL_i)$ acting with finite orbits on $\LL$ must be compact in the {\bf uniform}
topology on $\partial G$. 

\smallskip

\noindent As an approach to this, we propose {\it an alternate topology} on $PPQI(G,H)$
and call it the {\bf permutation topology}. Enumerate
$\LL = L_1, L_2, \cdots$. Since the representation of $PPQI(G,H)$ in the symmetric group
of permutations $S(\LL )$ is faithful, we declare that a system of open neighborhoods of the origin in 
$PPQI(G,H)$
is given by the set $U_N$ of elements of $PPQI(G,H)$ fixing $L_i, i = 1 \cdots N$. Now consider an 
element $\phi \in PPQI(G,H)$ acting with finite orbits on $\LL$. Then the (closed) monothetic
subgroup $(\phi )$ generated by $\phi$ is locally compact and by Corollary \ref{hscor}, under certain
hypotheses, it cannot contain the $p$-adics. Hence it must have {\bf arbitrarily small torsion
elements}. We cannot apply Theorem
\ref{ns} right away. To be able to apply Theorem
\ref{ns}, we need to show the following:\\ {\it For $\epsilon$ as in Theorem
\ref{ns}, there exists $N$ such that for all $k$, if $\phi^k$ stabilizes each $L_i, i = 1 \cdots N$, then
each orbit of $\phi^k$ has diameter less than $\epsilon$. }\\
Thus a weak enough statement ensuring a comparison of the {\it permutation topology} 
with the {\bf uniform topology} is necessary. The coarse barycenter construction of \cite{mahan-relrig}
might be helpful here to construct quasi-isometries coarsely fixing large balls and providing a starting
point for the problem. 

\smallskip

\noindent {\bf Problem 2:}  A more important and difficult
problem is to rule out elements of $PPQI(G,H)$ which fix finitely many elements of $\LL$ (and hence coarsely
fix the origin in $G$ by Theorem \ref{origin})
 but act with at least one unbounded orbit on $\LL$. We would have to show 
the following:\\ {\it There exists $N$ such that if $\phi \in PPQI(G,H)$ stabilizes each $L_i, i = 1 \cdots N$, then
each orbit of $\phi$ is finite. }

\subsection{Acknowledgments} I would like to thank Mladen Bestvina, Kingshook Biswas,  Marc Bourdon,
Emmanuel Breuillard and George Willis for  helpful conversations and
 correspondence.  Special thanks are due to Misha Gromov
for an extremely instructive conversation during which he told me about Scholium \ref{gromov}. 

This work was done in part while the author was visiting Caroline Series
at the Mathematics Institute, University of Warwick; and the
Universite Paris-Sud, Orsay under the Indo-French
ARCUS programme. I thank these institutions for their hospitality.

\bibliography{wpattern}
\bibliographystyle{alpha}

\noindent School of Mathematical Sciences, RKM Vivekananda University, \\
Belur Math, WB 711202, India \\
\texttt{mahan@rkmvu.ac.in}; \texttt{mahan.mj@gmail.com}

\end{document}